# HEAVY TRAFFIC LIMIT FOR A PROCESSOR SHARING QUEUE WITH SOFT DEADLINES

By H. Christian Gromoll[1] and Łukasz Kruk[2]

*Stanford University and Maria Curie-Skłodowska University*

This paper considers a GI/GI/1 processor sharing queue in which jobs have soft deadlines. At each point in time, the collection of residual service times and deadlines is modeled using a random counting measure on the right half-plane. The limit of this measure valued process is obtained under diffusion scaling and heavy traffic conditions and is characterized as a deterministic function of the limiting queue length process. As special cases, one obtains diffusion approximations for the lead time profile and the profile of times in queue. One also obtains a snapshot principle for sojourn times.

**1. Introduction.** Congestion and delay is typical in heavily loaded queueing systems and may vary significantly according to the service discipline being used. There is much interest in identifying service disciplines that minimize delay and this has led to a rich literature on performance analysis. Classical measures of delay include workload, queue length and sojourn time. In many circumstances, these performance measures help assess a service discipline's effectiveness at minimizing the time jobs spend in the system (their sojourn times). However, these measures may be inadequate for systems with more specific timing requirements. For example, in systems which have heterogeneous job deadlines, minimizing sojourn times may not be the best strategy for minimizing missed deadlines. Recent work has investigated various queueing models that allow for a more general notion of delay [4, 6, 13, 14, 15, 16, 18, 26]. In these models, jobs have individual

Received July 2005; revised June 2006.

[1]Supported in part by an NSF Mathematical Sciences Postdoctoral Research Fellowship, a European Union Marie Curie Postdoctoral Research Fellowship and EURANDOM.

[2]Supported by the State Committee for Scientific Research of Poland Grant 2 PO3A 012 23.

random deadlines and a service discipline's performance depends on how effectively it meets them.

Such *real-time* queueing models arise naturally in applications such as manufacturing systems, voice and video communication systems and avionics and automotive control systems. The associated performance measures are high-dimensional because individual timing information must be tracked to determine if deadlines are being met. Exact analysis is usually intractable, but heavy traffic approximations often exist. This was verified for the EDF (earliest deadline first) service discipline in [6] for the GI/GI/1 queue, in [15] for the multiclass queue with HLPS (head-of-the-line processor sharing) across classes, in [26] for feed-forward networks and in [16] for multiclass acyclic networks. Corresponding results for the FIFO (first in first out) discipline were also given in [15, 26]. The accuracy of these approximations in the single server case was investigated in [13, 14].

Processor sharing is a widely used idealization of the round-robin and time-sharing protocols used in computer and communication systems. Due to its practical importance, it is natural to investigate analogues of the results in [6, 15, 16, 26] for this service discipline. This is the subject of the present paper, which provides a framework for evaluating the heavy traffic performance of the GI/GI/1 processor sharing queue with respect to a fairly general structure of deadlines.

A new feature of this framework is that job deadlines are allowed to be correlated with service times; previous work has assumed independence of the two. This makes it possible, for example, to model scenarios in which large jobs have longer deadlines than small jobs. Although this adds a level of realism, it also requires some additional technical machinery. The state of the processor sharing queue with deadlines will be tracked using a measure valued process in the right half-plane. This idea builds on previous work on the GI/GI/1 processor sharing queue [8, 9, 10, 20, 21]. The setup requires detailed information about how arriving jobs affect the system state. Consequently, the primitive processes must be considered jointly as a measure valued arrival process. This creates some difficulty. The usual estimates needed to control scaling limits of an arrival process derive from functional weak laws of large numbers and functional central limit theorems. Analogous estimates for the arrival process used here require stronger results: the paper will draw on the theory of empirical processes to obtain functional Glivenko–Cantelli-type estimates for the measure valued arrival process (see Section 4.2).

Although dealing with measure valued processes in the half-plane requires some theoretical overhead, the additional work pays off. Section 3 provides some examples of the types of computation the more general setup allows. In particular, Section 3.2 uses it to compute a processor sharing analogue of Reiman's snapshot principle [18, 22] and Section 3.3 gives an example of how one might set realistic deadlines under the assumption of linear dependence between service times and deadlines. Furthermore, the present analysis for processor sharing is just an example of how empirical process theory may be applied to measure valued state descriptors. It is likely that the approach can be adapted to the study of other models.

The model considered here consists of an infinite capacity buffer, to which jobs arrive according to a delayed renewal process. Jobs arrive with service time requirements determined by a sequence of independent, identically distributed positive random variables. A single server processes buffered jobs at unit rate according to the processor sharing discipline: it works simultaneously on all jobs in the buffer, providing an equal fraction of its capacity to each. That is, if there are $Z(t) \geq 1$ jobs in the buffer at time $t$, then each job is receiving service at the instantaneous rate $1/Z(t)$. When the server has fulfilled a given job's service time requirement, the job exits the system. In addition to its service time, each job arrives at the buffer with a deadline given by the job's arrival time plus a real valued random variable called the *initial lead time*. The sequence of initial lead times is independent and identically distributed, but a job's initial lead time may be correlated with its service time. Job deadlines are *soft*, meaning that jobs remain in the system until served to completion; a job still in the system when its deadline expires is called *late*.

Since deadlines are soft and the server does not take them into account, the processor sharing discipline is unaltered by the inclusion of timing information in the model. In particular, the workload and queue length processes are identical to those of the classical processor sharing model. Thus, the result described below should be interpreted as a performance analysis of classical processor sharing, where the performance is measured with respect to the aforementioned deadline structure. In particular, the result answers the following question. Given a certain structure of deadlines, how well does the processor sharing discipline perform in meeting them? The setup used allows a variety of functionals to be considered in answering this question. For example, the lead time profile, time-in-queue profile and empirical sojourn time distributions can all be described using the main result (see Section 3). In this regard, the present analysis may be used to compare the effectiveness of processor sharing with that of other disciplines under the same deadline structure. As mentioned above, performance analyses of this type exist for certain FIFO models (another discipline that ignores the deadlines) as well as certain EDF models (a discipline which takes them into account).

Note that the system studied here is different from a GI/GI/1 processor sharing queue with *firm* deadlines, in which jobs exit early if their deadline expires. The latter system exhibits very different behavior and is the subject of forthcoming work. For more on processor sharing queues, see [25] for a survey up to 1987 and [10] for a discussion of more recent work.

The main result of this paper concerns a measure valued process that tracks the state of the processor sharing queue with soft deadlines. At time $t \geq 0$, this state includes both the *residual service time* and the *lead time* of each buffered job: if $\mathrm{I}(t)$ denotes the set of jobs in the buffer at time $t$, then the residual service time $v_i(t)$ is the remaining service requirement of job $i \in \mathrm{I}(t)$ and the lead time $l_i(t)$ is the deadline of job $i$ minus the current time $t$. If $l_i(t)$ is non-negative, it represents the remaining time until the deadline of job $i$ expires; if it is negative, its absolute value is the time overdue. Let $\mathbf{M}$ be the set of finite, non-negative Borel measures on the right half-plane $\mathbb{H}_+ = [0, \infty) \times (-\infty, \infty)$ and let $\delta_{(x,y)}$ denote the Dirac point measure at $(x, y) \in \mathbb{H}_+$. For each $t \geq 0$, the state descriptor $\mathcal{Z}(t)$ is the random element of $\mathbf{M}$ given by

$$\mathcal{Z}(t) = \sum_{i \in \mathrm{I}(t)} \delta_{(v_i(t), l_i(t))}.$$

Note that projection of $\mathcal{Z}(\cdot)$ onto the first coordinate yields the state descriptor used in [8, 9, 10, 20, 21]. For each $r$ in a sequence $\mathcal{R}$ of positive real numbers tending to infinity, define a diffusion scaled version of the state descriptor by

$$\hat{\mathcal{Z}}^r(t)(B \times C) = \frac{1}{r} \mathcal{Z}(r^2 t)(B \times rC), \qquad B \subset [0, \infty), \ C \subset (-\infty, \infty).$$

Let $\alpha$ be the limiting arrival rate (as $r \to \infty$) for jobs entering the system and let $\vartheta \in \mathbf{M}$ be the limiting joint distribution (as $r \to \infty$) of service times and initial lead times (the precise form of these assumptions is given in Section 2.3 below). For $z > 0$, let $\vartheta_e^z \in \mathbf{M}$ be the measure defined by

$$\vartheta_e^z([x, \infty) \times [y, \infty)) = \alpha \int_0^\infty \vartheta([x + uz^{-1}, \infty) \times [y + u, \infty)) \, du$$

for all $(x, y) \in \mathbb{H}_+$ and let $\vartheta_e^0 = \mathbf{0}$, where $\mathbf{0}$ denotes the zero measure in $\mathbf{M}$. The main result of this paper states that under mild conditions, including standard heavy traffic assumptions, the sequence of diffusion scaled state descriptors $\{\hat{\mathcal{Z}}^r(\cdot)\}$ converges in distribution as $r \to \infty$ to the measure valued process $\vartheta_e^{Z^*(\cdot)}$, where $Z^*(\cdot)$ is a reflected Brownian motion on $[0, \infty)$. The process $Z^*(\cdot)$ is the weak limit of the diffusion scaled queue length process obtained in [9]. Note that for each $z > 0$, the measure $\vartheta_e^z$ can be computed in a straightforward way, either numerically, via simulations, or, in some cases, explicitly. Thus, the process $\vartheta_e^{Z^*(\cdot)}$ is a tractable approximation to the dynamics of the processor sharing queue with soft deadlines.

The proof of this result proceeds in several stages and is summarized informally as follows. For each $t \geq 0$, define a fluid scaled version of the state descriptor by

$$\bar{\mathcal{Z}}^r(t)(B \times C) = \frac{1}{r} \mathcal{Z}(rt)(B \times rC), \qquad B \subset [0, \infty), \ C \subset (-\infty, \infty),$$

and for $m = 0, 1, 2, \ldots$, define shifted processes $\bar{\mathcal{Z}}^{r,m}(t) = \bar{\mathcal{Z}}^r(m + t)$. Fix times $T, L > 1$. It is first proved that $\mathcal{R}$-indexed sequences of sample paths of the shifted processes $\{\bar{\mathcal{Z}}^{r,m}(t) : r \in \mathcal{R}, m \le \lfloor rT \rfloor, t \in [0, L)\}$ are precompact with high probability. Next, limit points of these sequences are characterized as *local fluid limits*, which are continuous measure valued functions $\zeta : [0, L) \to \mathbf{M}$. In particular, it is shown that these local fluid limits uniformly approximate the sample paths of $\bar{\mathcal{Z}}^{r,m}(\cdot)$ on $[0, L)$ with asymptotically high probability. Moreover, for sufficiently large $L$ and $t \in [L-1, L)$, every local fluid limit $\zeta(\cdot)$ is approximately in steady state, which is given by

$$\zeta(t) \approx \vartheta_e^{z(t)}, \tag{1.1}$$

where $z(t) = \zeta(t)(\mathbb{H}_+)$ is the total mass of $\zeta(t)$.

For each $r \in \mathcal{R}$, the interval $[0, rT]$ is covered by the overlapping time intervals $[m, m + L)$, where $m = 0, \ldots, \lfloor rT \rfloor$. Thus, for any $t \in [0, rT]$, there exist $m \le \lfloor rT \rfloor$ and $s \in [0, L)$ such that $\bar{\mathcal{Z}}^r(t) = \bar{\mathcal{Z}}^{r,m}(s)$. Since $\hat{\mathcal{Z}}^r(t) = \bar{\mathcal{Z}}^r(rt)$, the sample paths of $\{\bar{\mathcal{Z}}^{r,m}(\cdot) : m \le \lfloor rT \rfloor\}$ on the time interval $[0, L)$ determine the sample paths of $\hat{\mathcal{Z}}^r(\cdot)$ on $[0, T]$. Consequently, (1.1) translates to the sample paths of $\hat{\mathcal{Z}}^r(\cdot)$ on $[(L-1)r^{-1}, T]$. That is, as $r \to \infty$,

$$\hat{\mathcal{Z}}^r(t) \approx \vartheta_e^{\hat{Z}^r(t)}, \qquad t \in [(L-1)r^{-1}, T], \tag{1.2}$$

where $\hat{Z}^r(\cdot) = \hat{\mathcal{Z}}^r(\cdot)(\mathbb{H}_+)$ is the diffusion scaled queue length process. With some extra work, one obtains (1.2) for all $t \in [0, T]$.

In particular, the diffusion scaled state descriptor $\hat{\mathcal{Z}}^r(\cdot)$ can be asymptotically recovered from the one-dimensional diffusion scaled queue length process $\hat{Z}^r(\cdot)$ by the lifting map $\Delta_\vartheta : z \mapsto \vartheta_e^z$. This phenomenon is known as *state space collapse*. The process $\hat{Z}^r(\cdot)$ converges in distribution to a reflected Brownian motion $Z^*(\cdot)$ by Corollary 2.4 of [9]. Applying the continuous mapping theorem to (1.2) and $\Delta_\vartheta$ completes the proof.

The program outlined above for proving state space collapse of a diffusion scaled process using overlapping sections of a fluid scaled process is motivated by methodology developed in [3, 24] for open multiclass networks with HL (head-of-the-line) service disciplines. This methodology was adapted to the measure valued state descriptor of a processor sharing queue in [9] and has been further adapted here. The limiting diffusion process $\Delta_\vartheta Z^*(\cdot) = \vartheta_e^{Z^*(\cdot)}$ takes values in a submanifold $\{\vartheta_e^z : z \ge 0\}$ of $\mathbf{M}$, known as the *invariant manifold*. Note that this manifold is not a linear space. This is an interesting contrast to other results exhibiting state space collapse, including [3, 5, 9, 24], in which the invariant manifold is frequently a linear space.

The paper is organized as follows. Section 2 gives a precise description of the model, different forms of scaling and various asymptotic assumptions;

it also contains the statement of the main theorem. Section 3 illustrates several applications of the theorem. Section 4 contains the proof of tightness of the shifted fluid scaled state descriptors. Section 5 investigates properties of local fluid limits and Section 6 contains the proof of state space collapse, leading quickly to the proof of the main theorem.

1.1. *Notation.* The following notation will be used throughout. Let $\mathbb{N} = \{1, 2, \ldots\}$ and let $\mathbb{Z}$ and $\mathbb{R}$ denote the integers and real numbers, respectively. Let $\mathbb{R}_+ = [0, \infty)$ and denote the right half-plane $\mathbb{R}_+ \times \mathbb{R}$ by $\mathbb{H}_+$. For $a, b \in \mathbb{R}$, write $a \vee b$ for the maximum, $a \wedge b$ for the minimum, $a^+$ and $a^-$ for the positive and negative parts, $\lfloor a \rfloor$ for the integer part and $\lceil a \rceil$ for the smallest integer $n \geq a$. Denote the indicator of a set $B \subset \mathbb{H}_+$ by $1_B$ and let $B^\varepsilon = \{w \in \mathbb{H}_+ : \inf_{z \in B} \|w - z\| < \varepsilon\}$. For $w \in \mathbb{H}_+$, let $B + w = \{z + w : z \in B\}$. Denote by $\mathbf{C}_b(S)$ the space of continuous bounded real valued functions on a topological space $S$.

Recall that $\mathbf{M}$ is the set of nonnegative finite Borel measures on $\mathbb{H}_+$, with zero measure denoted by $\mathbf{0}$. For $\zeta \in \mathbf{M}$ and a $\zeta$-integrable function $g : \mathbb{H}_+ \to \mathbb{R}$, define $\langle g, \zeta \rangle = \int_{\mathbb{H}_+} g \, d\zeta$. The space $\mathbf{M}$ is endowed with the weak topology: $\zeta_n \xrightarrow{\mathbf{w}} \zeta$ in $\mathbf{M}$ if and only if $\langle g, \zeta_n \rangle \to \langle g, \zeta \rangle$ for all $g \in \mathbf{C}_b(\mathbb{H}_+)$. With this topology, $\mathbf{M}$ is a Polish space [19]. It will be convenient to use the following metric: for $\xi, \zeta \in \mathbf{M}$, let

$$
\begin{aligned}
\mathbf{d}[\xi, \zeta] = \inf\{\varepsilon > 0 : &\, \xi(B) \leq \zeta(B^\varepsilon) + \varepsilon \text{ and} \\
&\, \zeta(B) \leq \xi(B^\varepsilon) + \varepsilon, \text{ for all closed } B \subset \mathbb{H}_+\}.
\end{aligned}
$$
(1.3)

It is straightforward to verify that $\mathbf{d}[\cdot, \cdot]$ is a complete metric on $\mathbf{M}$ inducing the weak topology.

All stochastic processes are assumed to be right continuous with finite left limits. Let $\mathbf{D}([0, \infty), \mathbf{M})$ denote the space of right continuous, left limited paths $\zeta : [0, \infty) \to \mathbf{M}$, endowed with the Skorohod $J_1$-topology. For two random objects $X, Y$ with the same distribution $\mu$, write $X \sim Y \sim \mu$. Write $X_n \Rightarrow Y$ when $X_n$ converges in distribution to $Y$ as $n \to \infty$.

**2. A processor sharing queue with soft deadlines.** This section gives a precise description of the model under consideration, specifies assumptions and states the main result. A formal definition of the processor sharing queue was given in [10] and the reader is referred there for a more detailed description. The present model is a generalization of the one considered in [10], so parts of the definition are restated here. This is the subject of Section 2.1 below, which introduces a sequence of processor sharing models incorporating soft deadlines, along with the associated notation. Section 2.2 describes the scaling and time shifts to be applied to this sequence and Section 2.3 specifies asymptotic assumptions. The main result appears in Section 2.4.

2.1. *Sequence of models.* Let $\mathcal{R} \subset (0, \infty)$ be a sequence that increases to infinity. Suppose that, for each $r \in \mathcal{R}$, there is a stochastic model consisting of the following: a processor sharing server that processes jobs at unit rate from an infinite capacity buffer, a collection of stochastic primitives $(E^r(\cdot), \{v_i^r, l_i^r\}_{i=1}^\infty)$ describing the arrival times, service times and initial lead times of jobs arriving to the buffer and a random initial condition specifying the state of the system at time 0. The time evolution of the system state is described by a collection of performance processes, defined in terms of the primitives and initial condition through a set of descriptive equations. The random objects in each model are defined on a probability space $(\Omega^r, \mathscr{F}^r, \mathbf{P}^r)$, with expectation on this space denoted $\mathbf{E}^r$.

The stochastic primitives consist of an *exogenous arrival process* $E^r(\cdot)$ and a sequence of *service times* and *initial lead times* $\{v_i^r, l_i^r\}_{i=1}^\infty$. The arrival process $E^r(\cdot)$ is a rate $\alpha^r$ delayed renewal process associated with a sequence $\{u_i^r\}_{i=1}^\infty$ of finite nonnegative *interarrival times*. For $t \geq 0$, $E^r(t)$ represents the number of jobs that have arrived at the buffer during the time interval $(0, t]$. The $i$th job to arrive after time 0 is called job $i$; jobs already in the buffer at time 0 will be called *initial jobs*. The quantity $u_1^r$ is the arrival time of the first job and for $i \geq 2$, $u_i^r$ is the elapsed time between the arrival of job $i - 1$ and job $i$. Thus, job $i$ arrives at time $U_i^r = \sum_{j=1}^i u_j^r$ for $i \geq 1$. Define $U_0^r = 0$. Then for $t \geq 0$,

$$E^r(t) = \sup\{i \geq 0 : U_i^r \leq t\}.$$

Assume that $\{u_i^r\}_{i=1}^\infty$ is a sequence of independent random variables and that $\{u_i^r\}_{i=2}^\infty$ are independent and identically distributed with mean $1/\alpha^r \in (0, \infty)$ and standard deviation $a^r < \infty$. The first element $u_1^r$ of the sequence is assumed to be strictly positive with finite mean.

For each $i \geq 1$, the service time $v_i^r$ represents the amount of processing time that job $i$ requires from the server; the initial lead time $l_i^r$ represents the maximum amount of time that job $i$ can be in the buffer without being late. Since this job arrives at time $U_i^r$, it will be late after time $U_i^r + l_i^r$. Late jobs remain in the buffer until completing service. Consequently, $U_i^r + l_i^r$ is a soft deadline for job $i$. Assume that $\{v_i^r\}$ are strictly positive random variables and that $\{v_i^r, l_i^r\}_{i=1}^\infty$ is a sequence of independent and identically distributed random vectors with common joint distribution given by a Borel probability measure $\vartheta^r$ on $\mathbb{H}_+$. Assume that $v_1^r$ has mean $1/\beta^r \in (0, \infty)$ and standard deviation $b^r < \infty$ and that $\{v_i^r, l_i^r\}_{i=1}^\infty$ is independent of the arrival process $E^r(\cdot)$.

It will be convenient to express the primitives as a measure valued arrival process. For each $t \geq 0$, let $\mathcal{L}^r(t) \in \mathbf{M}$ be the random measure

$$\mathcal{L}^r(t) = \sum_{i=1}^{E^r(t)} \delta_{(v_i^r, l_i^r)},$$

where $\delta_{(x,y)}$ denotes the Dirac measure at $(x,y) \in \mathbb{H}_+$.

The initial condition specifies $Z^r(0)$, the number of initial jobs present in the buffer at time zero, as well as the service time requirements and initial lead times of these initial jobs. Assume that $Z^r(0)$ is a nonnegative, integer valued random variable. The service times and initial lead times for initial jobs are the first $Z^r(0)$ elements of a sequence $\{\tilde{v}_j^r, \tilde{l}_j^r\}_{j=1}^\infty$ of random vectors, where $\{\tilde{v}_j^r\}$ are strictly positive. Assume that $Z^r(0)$ and $\{\tilde{v}_j^r, \tilde{l}_j^r\}_{j=1}^\infty$ are independent of $\{u_i^r\}_{i=2}^\infty$ and $\{v_i^r, l_i^r\}_{i=1}^\infty$; they do not need to be independent of each other. A convenient way to express the initial condition is to define an initial random measure $\mathcal{Z}^r(0) \in \mathbf{M}$ by

$$\mathcal{Z}^r(0) = \sum_{j=1}^{Z^r(0)} \delta_{(\tilde{v}_j^r, \tilde{l}_j^r)}.$$

Henceforth, $\mathcal{Z}^r(0)$ will be used as the initial condition. Let $\chi \colon \mathbb{H}_+ \to \mathbb{R}_+$ denote the projection $(x,y) \mapsto x$. Assume that $\mathcal{Z}^r(0)$ satisfies

$$(2.1) \qquad\qquad \mathbf{E}^r[\langle 1, \mathcal{Z}^r(0)\rangle] < \infty,$$

$$(2.2) \qquad\qquad \mathbf{E}^r[\langle \chi, \mathcal{Z}^r(0)\rangle] < \infty.$$

Note that since $\langle 1, \mathcal{Z}^r(0)\rangle = Z^r(0)$ and $\langle \chi, \mathcal{Z}^r(0)\rangle = \sum_{j=1}^{Z^r(0)} \tilde{v}_j^r$, assumptions (2.1) and (2.2) mean that the expected initial queue length and expected initial workload are finite.

The performance processes describe the system's behavior; they are defined in terms of the primitives through a set of equations that embody the processor sharing discipline and the expiration of deadlines. For each $t \geq 0$, let $S^r(t)$ denote the *cumulative service per job* provided by the server up to time $t$. Thus, if job $i \leq E^r(t)$ arrives at time $U_i^r \leq t$, the cumulative amount of processing time job $i$ receives by time $t$ is equal to $v_i^r \wedge (S^r(t) - S^r(U_i^r))$. Similarly, the cumulative amount of processing time an initial job $j \leq Z^r(0)$ receives by time $t$ is equal to $\tilde{v}_j^r \wedge S^r(t)$. Define the *residual service time* at time $t$ of job $i$ (and initial job $j$) by

$$(2.3) \qquad v_i^r(t) = (v_i^r - S^r(t) + S^r(U_i^r))^+, \qquad i = 1, \dots, E^r(t),$$

$$(2.4) \qquad \tilde{v}_j^r(t) = (\tilde{v}_j^r - S^r(t))^+, \qquad j = 1, \dots, Z^r(0).$$

A job's residual service time represents the current remaining amount of processing time required to fulfill its service time requirement. If a job's residual service time is zero, it has completed service and departed the buffer. Define the *lead time* at time $t$ of job $i$ (and of initial job $j$) by

$$(2.5) \qquad l_i^r(t) = U_i^r + l_i^r - t, \qquad i = 1, \dots, E^r(t),$$

$$(2.6) \qquad \tilde{l}_j^r(t) = \tilde{l}_j^r - t, \qquad j = 1, \dots, Z^r(0).$$

The lead time of a job currently in the buffer is the time remaining until its deadline passes; jobs with negative lead times are currently late.

The *state descriptor* $\mathcal{Z}^r(\cdot): [0, \infty) \to \mathbf{M}$ is a measure valued performance process that describes the time evolution of the system state: at each time $t \geq 0$, the random measure $\mathcal{Z}^r(t)$ has one unit of mass located at $(\tilde{v}_j^r(t), \tilde{l}_j^r(t)) \in \mathbb{H}_+$ for each initial job $j \leq Z^r(0)$ still in the buffer at time $t$ and one unit of mass at $(v_i^r(t), l_i^r(t)) \in \mathbb{H}_+$ for each job $i \leq E^r(t)$ still in the buffer at time $t$. For $(x, y) \in \mathbb{H}_+$, let $\delta_{(x,y)}^+$ be the Dirac measure at $(x, y)$ if $x > 0$, with $\delta_{(0,y)}^+ = \mathbf{0}$. Since a job is still in the buffer if and only if its current residual service time is positive, we have

$$(2.7) \qquad \mathcal{Z}^r(t) = \sum_{j=1}^{Z^r(0)} \delta_{(\tilde{v}_j^r(t), \tilde{l}_j^r(t))}^+ + \sum_{i=1}^{E^r(t)} \delta_{(v_i^r(t), l_i^r(t))}^+.$$

Let $Z^r(t)$ denote the number of jobs in the buffer, or queue length, at time $t \geq 0$. Then $Z^r(t) = \langle 1, \mathcal{Z}^r(t) \rangle$. Thus, under the processor sharing discipline, any job present in the buffer at time $t$ receives service at the instantaneous rate $\langle 1, \mathcal{Z}^r(t) \rangle^{-1}$. Note that if a job is present in the buffer at time $t$, then $\langle 1, \mathcal{Z}^r(t) \rangle \neq 0$. Let $\varphi(x) = 1/x$ for $x \in (0, \infty)$, with $\varphi(0) = 0$. Then the cumulative service per job up to time $t$ can be written

$$(2.8) \qquad S^r(t) = \int_0^t \varphi(\langle 1, \mathcal{Z}^r(s) \rangle) \, ds.$$

It will be convenient to let $S_{s,t}^r = S^r(t) - S^r(s)$ for $t \geq s \geq 0$.

Given the primitives $(E^r(\cdot), \{v_i^r, l_i^r\}_{i=1}^\infty)$ and the initial condition $\mathcal{Z}^r(0)$, the equations (2.3)–(2.8) determine the residual service times, lead times, the cumulative service per job process $S^r(\cdot)$ and the state descriptor $\mathcal{Z}^r(\cdot)$. This fact is not difficult, although somewhat tedious, to show.

Let $W^r(t)$ denote the workload in the buffer at time $t \geq 0$. This is the amount of time the server would have to work to complete the remaining service time requirements of all jobs in the buffer at time $t$, assuming no new arrivals take place. Since this equals the sum of the residual service times of all jobs present in the buffer at time $t$, we have

$$W^r(t) = \langle \chi, \mathcal{Z}^r(t) \rangle.$$

Note that the sequence of processor sharing models introduced here is a generalization of the sequence of models studied in [9, 10]. In particular, the measure valued process $\mu^r(\cdot) = \mathcal{Z}^r(\cdot) \circ \chi^{-1}$ is determined by the primitives $(E^r(\cdot), \{v_i^r\}_{i=1}^\infty)$ and the initial condition $(Z^r(0), \{\tilde{v}_j^r\}_{j=1}^\infty)$ and is the state descriptor that was studied in [9, 10].

2.2. *Scaling.* This paper concerns the asymptotic behavior, under heavy traffic conditions, of the $\mathcal{R}$-indexed sequence of models introduced in Section 2.1. To obtain useful limits as $r \to \infty$, various objects comprising the $r$th model must be appropriately scaled.

In Section 2.3 below, it is assumed that, as $r \to \infty$, the processor sharing models become heavily loaded at a rate governed by $r$. This implies that a job in the $r$th model remains in the buffer for a time that is of the order $r$ multiplied by the job's service time requirement. Thus, if initial lead times were to remain of the same order as service times, lead times $l_i^r(t)$, for large $t$ (of order $r^2$), would tend to $-\infty$ as $r \to \infty$. To obtain nontrivial scaling limits for the lead times, initial lead times $\{l_i^r\}$ in the $r$th model are assumed to be of order $r$ and will be scaled, along with lead times $l_i^r(\cdot)$, by $r^{-1}$. For each $r \in \mathcal{R}$, let $\breve{\vartheta}^r \in \mathbf{M}$ be the probability measure satisfying

$$(2.9) \qquad \breve{\vartheta}^r(B \times C) = \vartheta^r(B \times rC)$$

for all Borel sets $B \subset \mathbb{R}_+$ and $C \subset \mathbb{R}$.

The asymptotic heavy traffic behavior of the state descriptor will be examined on diffusion scale. For each $r \in \mathcal{R}$, the diffusion scaled state descriptor is defined, for $t \geq 0$, as the random measure $\hat{\mathcal{Z}}^r(t) \in \mathbf{M}$ satisfying

$$\hat{\mathcal{Z}}^r(t)(B \times C) = \frac{1}{r}\mathcal{Z}^r(r^2t)(B \times rC)$$

for all Borel sets $B \subset \mathbb{R}_+$ and $C \subset \mathbb{R}$. Note that this definition also scales lead times by $r^{-1}$.

Diffusion-scaled versions of the workload and queue length processes are also needed. For $t \geq 0$, define $\hat{W}^r(t) = r^{-1}W^r(r^2t)$ and $\hat{Z}^r(t) = r^{-1}Z^r(r^2t)$; note that $\hat{W}^r(t) = \langle \chi, \hat{\mathcal{Z}}^r(t) \rangle$ and $\hat{Z}^r(t) = \langle 1, \hat{\mathcal{Z}}^r(t) \rangle$. Since the workload process of a single server queue is the same for all work conserving service disciplines, including processor sharing, the heavy traffic behavior of $\hat{W}^r(\cdot)$ is described by the well-known result for FIFO queues [11]. Since $Z^r(\cdot) = \langle 1, \mathcal{Z}^r(\cdot) \rangle = \langle 1, \mathcal{Z}^r(\cdot) \circ \chi^{-1} \rangle$, the queue length process of the present model is the same as for the processor sharing queue without deadlines; a heavy traffic analysis of $\hat{Z}^r(\cdot)$ appears in [9]. The existing results for $\hat{W}^r(\cdot)$ and $\hat{Z}^r(\cdot)$ will be used in the analysis of $\hat{\mathcal{Z}}^r(\cdot)$.

Results for the diffusion scaled state descriptor $\hat{\mathcal{Z}}^r(\cdot)$ will be derived from results for the fluid scaled state descriptor, defined for $t \geq 0$ as the random measure $\bar{\mathcal{Z}}^r(t) \in \mathbf{M}$ satisfying

$$(2.10) \qquad \bar{\mathcal{Z}}^r(t)(B \times C) = \frac{1}{r}\mathcal{Z}^r(rt)(B \times rC)$$

for all Borel sets $B \subset \mathbb{R}_+$ and $C \subset \mathbb{R}$.

The relationship $\hat{\mathcal{Z}}^r(t) = \bar{\mathcal{Z}}^r(rt)$ will be essential to bootstrapping results from fluid scale up to diffusion scale. The process $\hat{\mathcal{Z}}^r(\cdot)$ is considered over a fixed finite time interval $[0, T]$ for $T > 1$. This corresponds to looking at the fluid scaled process $\bar{\mathcal{Z}}^r(\cdot)$ over $[0, rT]$. As $r \to \infty$, the asymptotic behavior of $\bar{\mathcal{Z}}^r(\cdot)$ over finite time intervals is analyzed using modified versions of the techniques developed in [9, 10]. Since $[0, rT]$ grows without bound as $r \to \infty$, it is necessary to piece together many (order $r$) overlapping sections of $\bar{\mathcal{Z}}^r(\cdot)$, each defined on a finite time interval of fixed length $L > 1$. This strategy is analogous to the one used in [3]. These overlapping sections of $\bar{\mathcal{Z}}^r(\cdot)$ are called the *shifted* fluid scaled state descriptors. For each $t \geq 0$ and $m \in \{0, 1, \dots\}$, define

$$(2.11) \qquad \bar{\mathcal{Z}}^{r,m}(t) = \bar{\mathcal{Z}}^r(m + t).$$

Then for each $r \in \mathcal{R}$, the time interval $[0, rT]$ is covered by the $\lfloor rT \rfloor + 1$ overlapping time intervals $[m, m + L]$ for $m = 0, \dots, \lfloor rT \rfloor$. Observe that for each $t \in [0, rT]$, there exists (at least one) $m \leq \lfloor rT \rfloor$ and $s \in [0, L]$ such that

$$\bar{\mathcal{Z}}^r(t) = \bar{\mathcal{Z}}^{r,m}(s).$$

Analysis of the processes $\bar{\mathcal{Z}}^{r,m}(\cdot)$ will involve fluid scaled and shifted fluid scaled versions of many of the processes introduced so far. For all $r \in \mathcal{R}$, $m \in \{0, 1, \dots\}$, $t \geq 0$ and $i = 1, \dots, E^r(rt)$, define

$$(2.12) \qquad \bar{E}^r(t) = \frac{1}{r} E^r(rt), \qquad \bar{E}^{r,m}(t) = \bar{E}^r(m + t),$$

$$(2.13) \qquad \bar{S}^r(t) = S^r(rt), \qquad \bar{S}^{r,m}(t) = \bar{S}^r(m + t),$$

$$\bar{Z}^r(t) = \frac{1}{r} Z^r(rt), \qquad \bar{Z}^{r,m}(t) = \bar{Z}^r(m + t),$$

$$(2.14) \qquad \bar{v}_i^r(t) = v_i^r(rt), \qquad \bar{v}_i^{r,m}(t) = \bar{v}_i^r(t + m),$$

$$(2.15) \qquad \bar{l}_i^r(t) = \frac{1}{r} l_i^r(rt), \qquad \bar{l}_i^{r,m}(t) = \bar{l}_i^r(t + m).$$

For $r \in \mathcal{R}$ and $t \geq 0$, define the fluid scaled measure valued arrival process by

$$(2.16) \qquad \bar{\mathcal{L}}^r(t) = \frac{1}{r} \sum_{i=1}^{r\bar{E}^r(t)} \delta_{(v_i^r, l_i^r r^{-1})}.$$

The fluid scaled processes $\bar{S}^r(\cdot)$ and $\bar{\mathcal{L}}^r(\cdot)$ will play particularly important roles; it will be convenient to have notation for their increments. Define the following fluid scaled and shifted fluid scaled increments: for all $r \in \mathcal{R}$, all $t \geq s \geq 0$ and all $m \in \{0, 1, \dots\}$, let

$$(2.17) \qquad \bar{S}_{s,t}^r = \bar{S}^r(t) - \bar{S}^r(s), \qquad \bar{S}_{s,t}^{r,m} = \bar{S}_{m+s,m+t}^r,$$

$$(2.18) \qquad \bar{\mathcal{L}}_{s,t}^r = \bar{\mathcal{L}}^r(t) - \bar{\mathcal{L}}^r(s), \qquad \bar{\mathcal{L}}_{s,t}^{r,m} = \bar{\mathcal{L}}_{m+s,m+t}^r.$$

Note that by (2.17) and (2.13),

$$\bar{S}_{s,t}^{r,m} = \int_{r(m+s)}^{r(m+t)} \varphi(\langle 1, \mathcal{Z}^r(u)\rangle)\, du = \int_s^t \varphi(\langle 1, \bar{\mathcal{Z}}^{r,m}(u)\rangle)\, du.$$

2.3. *Asymptotic assumptions.* This section imposes asymptotic assumptions on the sequence of models introduced in Section 2.1. This is the setting in which the main theorem is proved.

Let $\alpha$, $a$ and $p$ be fixed positive constants and let $\gamma \in \mathbb{R}$. Let $\vartheta$ be a probability measure on $\mathbb{H}_+$ satisfying

$$\vartheta(\{0\} \times \mathbb{R}) = 0, \tag{2.19}$$

$$\langle \chi^{4+p}, \vartheta \rangle < \infty, \tag{2.20}$$

$$\langle \chi, \vartheta \rangle = \alpha^{-1}. \tag{2.21}$$

Then $b = (\langle \chi^2, \vartheta \rangle - \langle \chi, \vartheta \rangle^2)^{1/2}$ is finite. For the sequence of arrival processes, assume that as $r \to \infty$,

$$(\alpha^r, a^r) \to (\alpha, a), \tag{2.22}$$

$$\mathbf{E}^r[u_1^r]/r \to 0, \tag{2.23}$$

$$\limsup_{r \to \infty} \mathbf{E}^r[(u_2^r)^{2+p}] < \infty. \tag{2.24}$$

For the sequence of service times and initial lead times, assume that as $r \to \infty$,

$$\breve{\vartheta}^r \xrightarrow{\mathbf{w}} \vartheta, \tag{2.25}$$

$$(\beta^r, b^r) \to (\alpha, b), \tag{2.26}$$

$$\limsup_{r \to \infty} \langle \chi^{4+p}, \breve{\vartheta}^r \rangle < \infty. \tag{2.27}$$

Define the *traffic intensity parameter* for the $r$th system by $\rho^r = \alpha^r/\beta^r$. Assumptions (2.22) and (2.26) imply that $\rho^r \to 1$, that is, the sequence of models approaches *heavy traffic*. Assume, further, that as $r \to \infty$,

$$r(1 - \rho^r) \to \gamma. \tag{2.28}$$

Assumption (2.28) represents the usual requirement that the sequence of systems approaches heavy traffic at a particular rate; the constant $\gamma$ appears as the drift coefficient of the limiting workload process. Assumption (2.23) implies that the initial residual interarrival time vanishes on diffusion scale. Assumption (2.25) specifies a weak limit for the joint distributions of service times and initial lead times when scaled according to (2.9). Assumptions (2.24) and (2.27) imply Lindeberg-type conditions which, along with (2.22), (2.23), (2.26) and (2.28), guarantee functional central limit theorems for the triangular arrays $\{u_i^r; i = 1, 2, \ldots\}_{r \in \mathcal{R}}$ and $\{v_i^r; i = 1, 2, \ldots\}_{r \in \mathcal{R}}$ and a functional Donsker theorem for the measure valued arrival processes $\mathcal{L}^r(\cdot)$. The central limit theorems ultimately imply that the diffusion scaled workload and queue length processes converge in distribution to reflected Brownian motions (see Propositions 4.1 and 4.2). The functional Donsker theorem implies a Glivenko–Cantelli-type estimate for $\mathcal{L}^r(\cdot)$. Assumption (2.27) is two moments stronger than is normally necessary for a functional central limit theorem. The additional restriction is required by Corollary 2.4 in [9], which will be applied in the present paper; see [9] for further discussion.

It remains to impose asymptotic assumptions on the diffusion scaled initial condition $\hat{\mathcal{Z}}^r(0)$. They are stated in terms of $\bar{\mathcal{Z}}^r(0) = \hat{\mathcal{Z}}^r(0)$ since they will be used in that form. The following definition is needed.

DEFINITION 2.1 (*Invariant manifold*). For each $z > 0$, let $\vartheta_e^z \in \mathbf{M}$ be the unique measure satisfying

$$(2.29) \qquad \langle 1_{[x,\infty) \times [y,\infty)}, \vartheta_e^z \rangle = \alpha \int_0^\infty \langle 1_{[x+uz^{-1},\infty) \times [y+u,\infty)}, \vartheta \rangle \, du$$

for all $x \in \mathbb{R}_+$ and $y \in \mathbb{R}$. Let $\vartheta_e^0 = \mathbf{0}$ and define

$$\mathbf{M}_\vartheta = \{\vartheta_e^z \in \mathbf{M} : z \geq 0\}.$$

By (2.21), (2.29) and a change of variables, $\langle 1, \vartheta_e^z \rangle = z$ for all $z \geq 0$. For $z > 0$, the probability measure $z^{-1}\vartheta_e^z$ can be thought of as the excess lifetime distribution of $\vartheta$ *in direction* $(z^{-1}, 1)$, a generalization of the notion of excess lifetime distribution to the half-plane $\mathbb{H}_+$. Following usage in [3], we call the one-parameter family of measures $\mathbf{M}_\vartheta \subset \mathbf{M}$ the *invariant manifold* associated with $\vartheta$. Note that the invariant manifold is not a linear space.

Let $\Theta \in \mathbf{M}$ be a random measure such that

$$(2.30) \qquad\qquad\qquad \Theta \in \mathbf{M}_\vartheta \qquad \text{a.s.},$$

$$(2.31) \qquad\qquad\qquad \mathbf{E}[\langle 1, \Theta \rangle] < \infty.$$

For the sequence of scaled initial measures $\bar{\mathcal{Z}}^r(0) = \hat{\mathcal{Z}}^r(0)$, assume that as $r \to \infty$,

$$(2.32) \qquad (\bar{\mathcal{Z}}^r(0), \langle \chi, \bar{\mathcal{Z}}^r(0) \rangle, \langle \chi^{1+p}, \bar{\mathcal{Z}}^r(0) \rangle) \Rightarrow (\Theta, \langle \chi, \Theta \rangle, \langle \chi^{1+p}, \Theta \rangle).$$

The second component of (2.32) implies that the scaled initial workload converges in distribution, that is, $\hat{W}^r(0) \Rightarrow \langle \chi, \Theta \rangle$.

The assumptions of this section are now summarized to simplify the statement of results that follow.

(A)
> *There is a sequence of processor sharing models, as
> defined in Section 2.1; there exist positive constants
> $\alpha$, $a$, $p$, a real number $\gamma \in \mathbb{R}$, a probability measure
> $\vartheta$ on $\mathbb{H}_+$ and a random measure $\Theta \in \mathbf{M}$ such that
> (2.19)–(2.32) hold.*

2.4. *Limit theorem.* Assume (A) and let $Z^*(\cdot)$ be a reflected Brownian
motion on $\mathbb{R}_+$ with drift $-2\gamma\alpha(1+\alpha^2 b^2)^{-1}$, variance $4\alpha^3(a^2+b^2)(1+\alpha^2 b^2)^{-2}$
and initial value $Z^*(0)$ that is equal in distribution to $\langle 1, \Theta \rangle$. Define the
measure valued process

$$\mathcal{Z}^*(\cdot) = \vartheta_e^{Z^*(\cdot)}.$$

THEOREM 2.2. *As $r \to \infty$, the sequence of diffusion scaled state descriptors $\{\hat{\mathcal{Z}}^r(\cdot)\}$ converges in distribution to the measure valued process $\mathcal{Z}^*(\cdot)$.*

Sections 4–6 are devoted to proving Theorem 2.2. Section 4 establishes
a tightness property for the families of shifted fluid scaled state descriptors
$\{\bar{\mathcal{Z}}^{r,m}(\cdot) : m \leq \lfloor rT \rfloor\}$. This yields the existence of limit points, called *local
fluid limits*, which are characterized in Section 5 as solutions of a certain
integral equation. In Section 6, the steady state behavior of these local fluid
limits, combined with the overlapping nature of the shifted fluid scaled state
descriptors $\{\bar{\mathcal{Z}}^{r,m}(\cdot) : m \leq \lfloor rT \rfloor\}$, reveals that $\{\hat{\mathcal{Z}}^r(\cdot)\}$ is asymptotically close
to the invariant manifold $\mathbf{M}_\vartheta$. Combined with an existing limit theorem
for the diffusion scaled queue length processes $\{\hat{Z}^r(\cdot)\}$, this establishes Theorem 2.2. The next section illustrates the applicability of Theorem 2.2 by
discussing several useful computations.

**3. Special cases and applications.** A central motivation of this paper is
to study the heavy traffic performance of the processor sharing discipline
in relation to timing requirements. A useful performance measure in this
context is called the *lead time profile* (see [6, 13, 14, 15, 16, 18, 26]). If
$\pi : \mathbb{H}_+ \to \mathbb{R}$ is the projection $(x, y) \mapsto y$, then the diffusion scaled lead time
profile is

$$\hat{\mathcal{Z}}^r(\cdot) \circ \pi^{-1}.$$

It is a random finite Borel measure on $\mathbb{R}$ that describes how well the service
discipline is meeting the soft deadlines of jobs. For example, the (scaled)
number of jobs that are currently late at time $t$ is $\langle 1_{(-\infty,0]}, \hat{\mathcal{Z}}^r(t) \circ \pi^{-1} \rangle$.

This section discusses a few special cases for which the lead time profile is
given by a convenient formula. When service times are independent of initial lead times, the lead time profile can be expressed in terms of a convolution of two measures on $\mathbb{R}$. If initial lead times equal zero, the lead time profile provides information about time in queue and sojourn times; see Sections 3.1 and 3.2. The case in which initial lead times depend linearly on service times leads to a heuristic rule for setting "realistic" deadlines for jobs in a processor sharing queue; see Section 3.3.

3.1. *Independence of service times and initial lead times.* Suppose that for each $r \in \mathcal{R}$, the service times and initial lead times of arriving jobs are independent, that is,

$$(3.1) \qquad \vartheta^r = \nu^r \times \lambda^r,$$

where $\nu^r$ is the service time distribution on $\mathbb{R}_+$ and $\lambda^r$ is the initial lead time distribution on $\mathbb{R}$. Letting $\breve{\lambda}^r(C) = \lambda^r(rC)$ for Borel sets $C \subset \mathbb{R}$, we have

$$\breve{\vartheta}^r = \nu^r \times \breve{\lambda}^r.$$

So, by (2.25),

$$\nu^r \times \breve{\lambda}^r \xrightarrow{\mathbf{w}} \vartheta \qquad \text{as } r \to \infty$$

and $\vartheta = \nu \times \lambda$, where $\nu$ and $\lambda$ are the weak limits of $\nu^r$ and $\breve{\lambda}^r$. By Theorem 2.2 and the continuous mapping theorem,

$$(3.2) \qquad \hat{\mathcal{Z}}^r(\cdot) \circ \pi^{-1} \Longrightarrow \vartheta_e^{Z^*(\cdot)} \circ \pi^{-1} \qquad \text{as } r \to \infty.$$

In the present setting, $\vartheta_e^0 = \mathbf{0}$ and for $z > 0$, $x \in \mathbb{R}_+$ and $y \in \mathbb{R}$,

$$(3.3) \qquad \langle 1_{(x,\infty)\times(y,\infty)}, \vartheta_e^z \rangle = \alpha \int_0^\infty \nu((x + uz^{-1}, \infty))\lambda((y + u, \infty))\, du.$$

To compute the limiting lead time profile, set $x = 0$ in (3.3) and observe that $\vartheta_e^z(\{0\} \times \mathbb{R}) = 0$. Then for $y \in \mathbb{R}$ and $z > 0$,

$$\langle 1_{\mathbb{R}_+\times(y,\infty)}, \vartheta_e^z \rangle = \alpha \int_0^\infty \nu((uz^{-1}, \infty))\lambda((y + u, \infty))\, du.$$

Thus, using Definition 2.1,

$$(3.4) \qquad \begin{aligned} \langle 1_{\mathbb{R}_+\times(-\infty,y]}, \vartheta_e^z \rangle &= \langle 1, \vartheta_e^z \rangle - \langle 1_{\mathbb{R}_+\times(y,\infty)}, \vartheta_e^z \rangle \\ &= z - \alpha \int_0^\infty \nu((uz^{-1}, \infty))\lambda((y + u, \infty))\, du. \end{aligned}$$

Since, by (2.21),

$$(3.5) \qquad z = z\alpha \int_0^\infty \nu((s, \infty))\, ds = \alpha \int_0^\infty \nu((uz^{-1}, \infty))\, du,$$

we can deduce from (3.4) that

$$(3.6) \qquad \langle 1_{\mathbb{R}_+ \times (-\infty, y]}, \vartheta_e^z \rangle = \alpha \int_0^\infty \nu((uz^{-1}, \infty)) \lambda((-\infty, y+u]) \, du.$$

Letting $\bar{\nu}_e^z$ be the Borel measure on $\mathbb{R}$ with density

$$(3.7) \qquad \bar{f}_e^z(u) = \begin{cases} \alpha \nu((-uz^{-1}, \infty)), & u \leq 0, \\ 0, & u > 0, \end{cases}$$

we can conclude from (3.6) that

$$(3.8) \qquad \langle 1_{\mathbb{R}_+ \times (-\infty, y]}, \vartheta_e^z \rangle = \int_{\mathbb{R}} \lambda((-\infty, y-u]) \, d\bar{\nu}_e^z(u).$$

This yields the convolution formula

$$(3.9) \qquad \vartheta_e^z \circ \pi^{-1} = \lambda \star \bar{\nu}_e^z.$$

Note that if $\nu_e$ denotes the excess lifetime distribution of $\nu$, then $\nu_e$ has density $\bar{f}_e^z(-uz)$ for $u \in \mathbb{R}$. That is, $\bar{\nu}_e^z$ is related to $\nu_e$ by a scaling factor $z$ and a reflection about zero.

Letting $\bar{\nu}_e^0$ be the zero measure on $\mathbb{R}$, (3.9) holds for all $z \geq 0$. Combining (3.2) and (3.9) yields an approximation result for the diffusion scaled lead time profiles. As $r \to \infty$,

$$(3.10) \qquad \hat{\mathcal{Z}}^r(\cdot) \circ \pi^{-1} \Longrightarrow \lambda \star \bar{\nu}_e^{Z^*(\cdot)}.$$

As a special case of (3.1), suppose, in addition, that the service time distribution $\nu^r$ is exponential for each $r \in \mathcal{R}$. Then $\nu$ is exponential with mean $\alpha^{-1}$. For $z > 0$, (3.7) takes the form

$$\bar{f}_e^z(u) = \begin{cases} \alpha \exp\left(\dfrac{\alpha}{z} u\right), & u \leq 0, \\ 0, & u > 0. \end{cases}$$

So, by (3.10), when $Z^*(t) > 0$, the limiting lead time profile $\lambda \star \bar{\nu}_e^{Z^*(t)}$ at time $t$ is equal to the convolution of $\lambda$ with an exponential measure with parameter $\alpha/Z^*(t)$ and total mass $Z^*(t)$ that is reflected about zero. This result applies, in particular, to a sequence of M/M/1 queues with lead times satisfying $\check{\lambda}^r = \lambda$ for all $r$. Thus, it verifies a conjecture of Lehoczky [18], who supported this statement with an informal analysis and Monte Carlo simulations.

3.2. *Time in queue, sojourn times and the snapshot.* An alternative special case of (3.1) can be used to study the profile of *times in queue* for the processor sharing discipline. A job's time in queue at time $t$ is the amount of time spent in the buffer by time $t$. On diffusion scale, the profile of times in queue $\hat{\tau}^r(t)$ at time $t$ is a random finite Borel measure on $\mathbb{R}_+$. It can be directly deduced from the diffusion scaled lead time profile by assuming that initial lead times are zero for all $r \in \mathcal{R}$. A job's time in queue at time $t$ then equals the absolute value of its current lead time; for $y \geq 0$,

$$(3.11) \qquad \langle 1_{[y,\infty)}, \hat{\tau}^r(\cdot) \rangle = \langle 1_{\mathbb{R}_+ \times (-\infty, -y]}, \hat{\mathcal{Z}}^r(\cdot) \rangle.$$

Since (3.1) holds when $\lambda^r = \delta_0$ for all $r \in \mathcal{R}$, the limiting lead time profile $\vartheta_e^{Z^*(\cdot)} \circ \pi^{-1}$ is given by (3.9) with $\lambda = \delta_0$; for $z = 0$, $\vartheta_e^z \circ \pi^{-1}$ equals the zero measure on $\mathbb{R}$ and for $z > 0$ and $y \in \mathbb{R}$, a short computation using (3.9) and (3.7) yields

$$(3.12) \qquad \langle 1_{\mathbb{R}_+ \times (-\infty, y]}, \vartheta_e^z \rangle = \begin{cases} z, & y \geq 0, \\ z\nu_e([-yz^{-1}, \infty)), & y < 0. \end{cases}$$

Thus, if $\tau^*(\cdot)$ denotes the limiting profile of times in queue, then for $t \geq 0$ and $y \geq 0$,

$$(3.13) \qquad \langle 1_{[y,\infty)}, \tau^*(t) \rangle = \begin{cases} Z^*(t)\nu_e([yZ^*(t)^{-1}, \infty)), & Z^*(t) > 0, \\ 0, & Z^*(t) = 0. \end{cases}$$

A related computation reveals the limiting sojourn time behavior, which is informally described below. A job's sojourn time is the time between its arrival at the buffer and its service completion. Consider the empirical distribution of a collection of sojourn times, sampled over a time interval that is small on diffusion scale. This distribution can be studied as the interval varies over diffusion scaled time. Since the interval will be taken to be asymptotically small, it can be defined in terms of residual service times.

Fix $T > 1$ and $t \in [0, T)$ such that $Z^*(t) > 0$. For $\varepsilon > 0$ and $r \in \mathcal{R}$, let $\hat{J}_\varepsilon^r(t)$ denote the set of jobs with residual service time in $(0, \varepsilon]$ at (diffusion scaled) time $t$. Let $\hat{\tau}_\varepsilon^r(t)$ be the diffusion scaled distribution of times in queue for those jobs in $\hat{J}_\varepsilon^r(t)$. Once again, assume that $\lambda^r = \delta_0$ for all $r$ so that a job's time in queue is the absolute value of its lead time. Then for $y \geq 0$,

$$\langle 1_{[0,y)}, \hat{\tau}_\varepsilon^r(t) \rangle = \varphi(\langle 1_{[0,\varepsilon] \times \mathbb{R}}, \hat{\mathcal{Z}}^r(t) \rangle) \langle 1_{[0,\varepsilon] \times (-y,\infty)}, \hat{\mathcal{Z}}^r(t) \rangle,$$

where $\varphi(x) = 1/x$ for $x > 0$, with $\varphi(0) = 0$. Let $\hat{\sigma}_\varepsilon^r(t)$ denote the diffusion scaled distribution of sojourn times for those jobs in $\hat{J}_\varepsilon^r(t)$. That is, if $\sigma_i^r$ is the unscaled sojourn time of job $i \in \hat{J}_\varepsilon^r(t)$, then for $y \geq 0$,

$$(3.14) \qquad \langle 1_{[0,y)}, \hat{\sigma}_\varepsilon^r(t) \rangle = \varphi(\langle 1_{[0,\varepsilon] \times \mathbb{R}}, \hat{\mathcal{Z}}^r(t) \rangle) \sum_{i \in \hat{J}_\varepsilon^r(t)} 1_{[0,y)}(r^{-1}\sigma_i^r).$$

For small $\varepsilon$, the probability measure $\hat{\tau}_\varepsilon^r(t)$ approximates $\hat{\sigma}_\varepsilon^r(t)$ in the following sense: if job $i \in \hat{J}_\varepsilon^r(t)$, then in the unscaled system, $v_i^r(r^2 t) \leq \varepsilon$. Since the sequence $\{\hat{\mathcal{Z}}^r(\cdot)\}$ is tight, there exists $C > 0$ such that, with high probability, $\sup_{s \in [0,T]} \hat{\mathcal{Z}}^r(s) \leq C$. Assume that $r$ is sufficiently large that $r^2 t + r\varepsilon C \leq r^2 T$. Then, with high probability, job $i$ receives service at rate at least $(rC)^{-1}$ if it is present in the buffer during the entire time interval $[r^2t, r^2t + r\varepsilon C] \subset [0, r^2T]$. It must therefore complete service at or before time $r^2t + r\varepsilon C$. Since its time in queue increases by at most $r\varepsilon C$ before it departs, its sojourn time $\sigma_i^r$ will differ from its time in queue at time $r^2t$ by at most $r\varepsilon C$. So, by (3.14), its diffusion scaled sojourn time will differ from its diffusion scaled time in queue at time $t$ by at most $\varepsilon C$. Taking $\varepsilon$ to be small compared to $C$, it is easy to show that $\hat{\tau}_\varepsilon^r(t)$ and $\hat{\sigma}_\varepsilon^r(t)$ are close in the weak topology. As $r \to \infty$, the number of jobs in $\hat{J}_\varepsilon^r(t)$ grows without bound for all $\varepsilon > 0$, yielding nontrivial distributions $\hat{\tau}_\varepsilon^r(t)$ and $\hat{\sigma}_\varepsilon^r(t)$.

By letting $\varepsilon \to 0$ after $r \to \infty$, one obtains a limiting distribution $\sigma^*(t)$. By Theorem 2.2 and the continuous mapping theorem, $\hat{\tau}_\varepsilon^r(t)$ can be approximated by $\tau_\varepsilon^*(t)$, where, for $y \geq 0$,

$$\langle 1_{[0,y)}, \tau_\varepsilon^*(t) \rangle = \langle 1_{[0,\varepsilon] \times (-y,\infty)}, \vartheta_e^{Z^*(t)} \rangle \varphi(\langle 1_{[0,\varepsilon] \times \mathbb{R}}, \vartheta_e^{Z^*(t)} \rangle).$$

For $z > 0$ and $y \geq 0$, a short computation using (3.3) with $\lambda = \delta_0$ yields

$$\frac{\langle 1_{[0,\varepsilon] \times (-y,\infty)}, \vartheta_e^z \rangle}{\langle 1_{[0,\varepsilon] \times \mathbb{R}}, \vartheta_e^z \rangle} = \frac{\nu_e([0, yz^{-1})) - \nu_e([\varepsilon, \varepsilon + yz^{-1}))}{\nu_e([0,\infty)) - \nu_e([\varepsilon,\infty))}$$

$$= \frac{\nu_e([0,\varepsilon)) - \nu_e([yz^{-1}, \varepsilon + yz^{-1}))}{\nu_e([0,\varepsilon))}.$$

Thus, for $z > 0$ and $y \geq 0$ such that $\nu(\{yz^{-1}\}) = 0$,

$$\lim_{\varepsilon \to 0} \frac{\langle 1_{[0,\varepsilon] \times (-y,\infty)}, \vartheta_e^z \rangle}{\langle 1_{[0,\varepsilon] \times \mathbb{R}}, \vartheta_e^z \rangle} = 1 - \nu([yz^{-1}, \infty)) = \nu([0, yz^{-1})).$$

Conclude that for $y \geq 0$,

$$\sigma^*(t)([0,y)) = \nu([0, yZ^*(t)^{-1})).$$

Since $\hat{\sigma}_\varepsilon^r(t)$ and $\hat{\tau}_\varepsilon^r(t)$ are close for small $\varepsilon$, the measure $\sigma^*(t)$ may be interpreted as a limiting sojourn time distribution of jobs "departing the system at diffusion scaled time $t$." Note that the order in which limits are taken here is the only way to obtain a nontrivial limit; if one takes $\varepsilon \to 0$ in the prelimit, then $\hat{\tau}_\varepsilon^r(t)$ converges to the zero measure because $\hat{\mathcal{Z}}^r(t)$ does not charge $\{0\} \times \mathbb{R}$. Also, note that the limit $\sigma^*(t)$ has been obtained for fixed $t$ when $Z^*(t) > 0$, not as a process.

An implication of (3.14) and the subsequent computations is that sojourn times in the $r$th system are of order $r$. Compare this to the state descriptor: since $\hat{\mathcal{Z}}^r(\cdot) = r^{-1}\mathcal{Z}^r(r^2\cdot)$ converges to a process with continuous sample paths, the state descriptor varies slowly on the time scale of sojourn times when $r$ is large. That is, jobs observe only very small changes in the system state, including the queue length, during their stay in the system. This is analogous to Reiman's *snapshot principle* (see [18, 22]).

3.3. *Linear dependence of initial lead times and service times.* By the results in the previous section, the amount of time a job spends in a heavily loaded single server processor sharing queue is roughly proportional to its service time and the queue length observed upon arrival. Suppose that in the $r$th system, jobs have "realistic" deadlines, that is, $l_i^r = crv_i^r$ for all $i = 1, 2, \ldots$, where $c > 0$ is a constant corresponding to the (diffusion scaled) queue length a job hopes to find upon arrival. Then for $x \in \mathbb{R}_+$ and $y \in \mathbb{R}$,

$$\langle 1_{[x,\infty) \times [y,\infty)}, \vartheta \rangle = \nu([x \vee yc^{-1}, \infty)).$$

By Theorem 2.2, $\hat{\mathcal{Z}}^r(\cdot) \Rightarrow \vartheta_e^{Z^*(\cdot)}$, where for $z \geq 0$, $x \in \mathbb{R}_+$ and $y \in \mathbb{R}$,

$$\langle 1_{[x,\infty) \times [y,\infty)}, \vartheta_e^z \rangle = z\alpha \int_0^\infty \langle 1_{[(x+u) \vee (y+zu)c^{-1}, \infty)}, \nu \rangle \, du.$$

Setting $x = 0$ yields, after some computation,

$$\langle 1_{\mathbb{R}_+ \times [y,\infty)}, \vartheta_e^z \rangle = \begin{cases} z, & y \leq 0, \\ c\nu_e([yc^{-1}, \infty)) - (c-z)\nu_e([y(c-z)^{-1}, \infty)), & y > 0, \end{cases}$$

for $0 < z \leq c$, where the last term above is understood to be zero when $z = c$. For $z > c$,

$$\langle 1_{\mathbb{R}_+ \times [y,\infty)}, \vartheta_e^z \rangle = \begin{cases} z + (c-z)\nu_e([y(c-z)^{-1}, \infty)), & y \leq 0, \\ c\nu_e([yc^{-1}, \infty)), & y > 0. \end{cases}$$

As the snapshot principle predicts, there is no lateness in the limiting system when $Z^*(t) \leq c$ and therefore almost no lateness in the prelimit when $r$ is large and $Z^r(t) \leq rc + o(r)$. If $Z^*(t) > c$ [$r$ is large and $Z^r(t) > r(c + \varepsilon)$ for some $\varepsilon > 0$], then jobs become late before they exit the system. Thus, $c$ acts as a threshold for $Z^*(\cdot)$, governing when the system exhibits lateness. Since the law of $Z^*(\cdot)$ is known [9], one can influence the proportion of missed deadlines by choosing $c$ appropriately. For example, in cases where $\gamma > 0$ (corresponding to a sequence of subcritical prelimit models), $Z^*(\cdot)$ has a steady state distribution [9]; by choosing $c$ as a particular quantile of this distribution, one could achieve a desired long-run average lateness rate.

**4. Tightness.** The remainder of the paper is devoted to proving Theorem 2.2. The goal of this section is to prove that $\mathcal{R}$-indexed sequences of sample paths chosen from $\{\bar{\mathcal{Z}}^{r,m}(\cdot) : m \leq \lfloor rT \rfloor\}$ are relatively compact with high probability.

4.1. *Preliminary results.* Assumption (A) has some important consequences which are reviewed below for later reference. Recall that $\hat{W}^r(\cdot) = \langle \chi, \hat{\mathcal{Z}}^r(\cdot) \rangle$ is the diffusion scaled workload process for the $r$th model. Let $W^*(\cdot)$ be a reflected Brownian motion on $\mathbb{R}_+$ with drift $-\gamma$, variance $\alpha(a^2 + b^2)$ and initial value $W^*(0)$ equal in distribution to $\langle \chi, \Theta \rangle$.

PROPOSITION 4.1.   *As $r \to \infty$, the sequence of diffusion scaled workload processes $\{\langle \chi, \hat{\mathcal{Z}}^r(\cdot) \rangle\}$ converges in distribution to $W^*(\cdot)$.*

PROOF.   It is well known that (A) is sufficient to imply the above result for the workload process of any single-server, single-class queue operating under a work conserving service discipline (including the present processor sharing model). The use of functional central limit theorems and continuous mappings to prove such results goes back to [11]. See [24] for a detailed account covering a more general setting than the present one.   □

A similar result holds for the diffusion scaled queue length process $\hat{Z}^r(\cdot) = \langle 1, \hat{\mathcal{Z}}^r(\cdot) \rangle$. Assume (A) and let $W^*(\cdot)$ be the reflected Brownian motion specified above. Let $C_\vartheta = 2\alpha(1 + \alpha^2 b^2)^{-1}$ and define $Z^*(\cdot) = C_\vartheta W^*(\cdot)$.

PROPOSITION 4.2.   *As $r \to \infty$, the sequence of diffusion scaled queue length processes $\{\hat{Z}^r(\cdot)\}$ converges in distribution to $Z^*(\cdot)$.*

PROOF.   See Theorem 2.3 in [9] and Corollary 2.4 in [9].   □

A simple corollary will also be needed. Assume (A) and let $T > 1$ and $\varepsilon, \eta \in (0, 1)$.

COROLLARY 4.3.   *For all $L > 1$,*

$$(4.1) \qquad \liminf_{r \to \infty} \mathbf{P}^r \bigg( \sup_{\substack{m \leq \lfloor rT \rfloor, \\ t \in [0, L]}} |\langle 1, \bar{\mathcal{Z}}^{r,m}(t) \rangle - \langle 1, \bar{\mathcal{Z}}^{r,m}(0) \rangle| \leq \varepsilon \bigg) \geq 1 - \eta.$$

*Furthermore, there exists $M > 1$ such that for all $L > 1$,*

$$(4.2) \qquad \liminf_{r \to \infty} \mathbf{P}^r \bigg( \sup_{\substack{m \leq \lfloor rT \rfloor, \\ t \in [0, L]}} \langle 1, \bar{\mathcal{Z}}^{r,m}(t) \rangle \vee \langle \chi, \bar{\mathcal{Z}}^{r,m}(t) \rangle \leq M \bigg) \geq 1 - \eta.$$

PROOF.   Expand the definition of $\bar{\mathcal{Z}}^{r,m}(\cdot)$ to rewrite in terms of diffusion scaling. It suffices to show that for all $L > 1$,

$$\liminf_{r \to \infty} \mathbf{P}^r \bigg( \sup_{\substack{t \in [0, T], \\ h \in [0, Lr^{-1}]}} |\hat{Z}^r(t + h) - \hat{Z}^r(t)| \leq \varepsilon \bigg) \geq 1 - \eta,$$

and that there exists $M > 1$ such that for all $L > 1$,

$$\liminf_{r \to \infty} \mathbf{P}^r \bigg( \sup_{t \in [0, T + Lr^{-1}]} \hat{Z}^r(t) \vee \langle \chi, \hat{\mathcal{Z}}^r(t) \rangle \leq M \bigg) \geq 1 - \eta.$$

Both statements follow from Propositions 4.1 and 4.2; the sequences $\{\langle \chi, \hat{\mathcal{Z}}^r(\cdot)\rangle\}$ and $\{\hat{Z}^r(\cdot)\}$ are tight and the limiting processes $W^*(\cdot)$ and $Z^*(\cdot)$ are a.s. continuous. Note that $M$ can indeed be chosen independent of $L$.  $\square$

The following result concerning the invariant manifold $\mathbf{M}_\vartheta$ will be needed.

DEFINITION 4.4 (*Lifting map*).   For a probability measure $\vartheta$ on $\mathbb{H}_+$ satisfying (2.19)–(2.21), let $\Delta_\vartheta : \mathbb{R}_+ \to \mathbf{M}_\vartheta$ be the associated lifting map given by

$$(4.3) \qquad\qquad \Delta_\vartheta z = \vartheta_e^z, \qquad z \in \mathbb{R}_+.$$

LEMMA 4.5.   *The map $\Delta_\vartheta : \mathbb{R}_+ \to \mathbf{M}_\vartheta$ is continuous and $\mathbf{M}_\vartheta \subset \mathbf{M}$ is closed.*

PROOF.   By (2.21), (2.29) and a change of variables, $\langle 1, \vartheta_e^y \rangle = y$ for all $y \geq 0$. Let $\{z_n\} \subset \mathbb{R}_+$ satisfy $z_n \to z$ as $n \to \infty$. Then $\langle 1, \vartheta_e^{z_n} \rangle \to \langle 1, \vartheta_e^z \rangle$. If $z = 0$, then $\vartheta_e^{z_n} \xrightarrow{\mathbf{w}} \mathbf{0}$. Assume without loss of generality that $z, z_n > 0$ for all $n$ and let $B \subset \mathbb{H}_+$ be closed. By (2.29) and two applications of Fatou's lemma,

$$\limsup_{n \to \infty} \vartheta_e^{z_n}(B) \leq \alpha \int_0^\infty \left\langle \limsup_{n \to \infty} 1_{B + (uz_n^{-1}, u)}, \vartheta \right\rangle du.$$

Since $B$ is closed,

$$\limsup_{n \to \infty} \vartheta_e^{z_n}(B) \leq \alpha \int_0^\infty \langle 1_{B + (uz^{-1}, u)}, \vartheta \rangle\, du = \vartheta_e^z(B).$$

Conclude from the Portmanteau theorem that $z_n^{-1} \vartheta_e^{z_n} \xrightarrow{\mathbf{w}} z^{-1} \vartheta_e^z$ and therefore, $\vartheta_e^{z_n} \xrightarrow{\mathbf{w}} \vartheta_e^z$. For the second statement, let $\{\vartheta_e^{z_n}\} \subset \mathbf{M}_\vartheta$ satisfy $\vartheta_e^{z_n} \xrightarrow{\mathbf{w}} \xi$ as $n \to \infty$. Then $z_n \to \langle 1, \xi \rangle$. Since $\Delta_\vartheta$ is continuous, $\xi = \vartheta_e^{\langle 1, \xi \rangle} \in \mathbf{M}_\vartheta$.  $\square$

4.2. *Functional Glivenko–Cantelli estimate.*   Proofs appearing below will frequently rely on a functional uniform Glivenko–Cantelli estimate for the measure valued arrival process specified in Section 2.1. They require uniform control over the empirical distributions of random vectors $(v_i^r, l_i^r r^{-1})$ associated with jobs arriving during an interval $[s, t]$. More explicitly, for all intervals $[s, t]$ with length bounded by a fixed constant $L$, the random measure

$$\bar{\mathcal{L}}_{s,t}^{r,m} = \frac{1}{r} \sum_{i = r\bar{E}^{r,m}(s)+1}^{r\bar{E}^{r,m}(t)} \delta_{(v_i^r, l_i^r r^{-1})}$$

must be asymptotically approximated, with high probability, by the measure $\alpha^r(t-s)\breve{\vartheta}^r$. For fixed $m$, $[s,t]$ and for a fixed Borel set $B \subset \mathbb{H}_+$, it follows from (A) and the weak law of large numbers that

$$(4.4) \qquad \bar{\mathcal{L}}_{s,t}^{r,m}(B) - \alpha^r(t-s)\breve{\vartheta}^r(B) \xrightarrow{\mathbf{P}^r} 0 \qquad \text{as } r \to \infty.$$

The present setting requires this approximation to be uniform over a sufficiently rich class $\mathscr{A}$ of subsets of $\mathbb{H}_+$ and, for each $T, L > 1$, uniform over all $m \le \lfloor rT \rfloor$ and all intervals $[s,t] \subset [0,L]$. This section addresses this issue using well-known results from empirical process theory.

Let $\mathscr{V}$ be a family of Borel-measurable functions $f : \mathbb{H}_+ \to \mathbb{R}$ and let $F : \mathbb{H}_+ \to \mathbb{R}$ be a Borel-measurable envelope for $\mathscr{V}$, that is, $|f| \le F$ for all $f \in \mathscr{V}$. For $f \in \mathscr{V}$ and a Borel probability measure $Q$ on $\mathbb{H}_+$, let $\|f\|_{Q,2} = (\int f^2\, dQ)^{1/2}$. Let $\ell^\infty(\mathscr{V})$ be the space of bounded functions $\mathcal{G} : \mathscr{V} \to \mathbb{R}$, equipped with the supremum norm $\|\mathcal{G}(f)\|_{\mathscr{V}} = \sup_{f \in \mathscr{V}} |\mathcal{G}(f)|$. For a random measure $\mathcal{L} \in \mathbf{M}$ (defined on $\Omega^r$) satisfying $\langle F, \mathcal{L} \rangle < \infty$ almost surely, identify $\mathcal{L}$ with a random element $\mathcal{L} \in \ell^\infty(\mathscr{V})$ satisfying $\mathcal{L}(f) = \langle f, \mathcal{L} \rangle$ almost surely. This identification is unique in law. Then for $r \in \mathcal{R}$, $m \le \lfloor rT \rfloor$ and $[s,t] \subset [0,L]$, the random measure $\bar{\mathcal{L}}_{s,t}^{r,m} - \alpha^r(t-s)\breve{\vartheta}^r$ also represents a random element of $\ell^\infty(\mathscr{V})$. A uniform version of (4.4) requires that for an appropriate family $\mathscr{V}$ of indicator functions, $\bar{\mathcal{L}}_{s,t}^{r,m} - \alpha^r(t-s)\breve{\vartheta}^r \xrightarrow{\mathbf{P}^r} \mathbf{0}$ in $\ell^\infty(\mathscr{V})$, uniformly in $m \le \lfloor rT \rfloor$ and $[s,t] \subset [0,L]$.

This uniform convergence holds under easily verifiable conditions on $\mathscr{V}$. The first main condition guarantees a minimal amount of measurability when dealing with suprema over $\mathscr{V}$. Call a family $\mathscr{V}$ a *Borel-measurable class* if for each $n \in \mathbb{N}$ and $(e_1, \ldots, e_n) \in \{-1, 1\}^n$, the map

$$(x_1, \ldots, x_n) \mapsto \left\| \sum_{i=1}^n e_i f(x_i) \right\|_{\mathscr{V}}$$

is Borel measurable on $\mathbb{H}_+^n$. The condition states that for all $\delta > 0$ and $r \in \mathcal{R}$, the families $\mathscr{V}_\delta^r = \{f - g : f, g \in \mathscr{V}, \|f - g\|_{\breve{\vartheta}^r, 2} < \delta\}$ and $\mathscr{V}_\infty^2 = \{(f-g)^2 : f, g \in \mathscr{V}\}$ are Borel-measurable classes.

The second main condition ensures that $\mathscr{V}$ is not too large, as measured by the following combinatorial property. Given a set $S$, a *Vapnik–Červonenkis class* (*VC-class*) is a collection $\mathscr{S}$ of subsets of $S$ with finite *Vapnik–Červonenkis index* (*VC-index*). The VC-index of a collection of subsets is the smallest integer $n$ ($\infty$ if none exists) such that $\mathscr{S}$ shatters no $n$-point subset $\{x_1, \ldots, x_n\} \subset S$. A collection $\mathscr{S}$ *shatters* $\{x_1, \ldots, x_n\}$ if each of the $2^n$ subsets of $\{x_1, \ldots, x_n\}$ can be *picked out* by $\mathscr{S}$, that is, can be written $A \cap \{x_1, \ldots, x_n\}$ for some $A \in \mathscr{S}$. Finally, a family of functions $\mathscr{U} : S \to \mathbb{R}$ is a VC-class if the collection of subgraphs $\{\{(x,y) : y < f(x)\} : f \in \mathscr{U}\}$ is a VC-class of sets in $S \times \mathbb{R}$. For a more detailed discussion of VC-classes, see, for example, [23], Chapter 2.6.

VC-classes satisfy the *uniform entropy condition*: Let $\mathscr{V}$ be a VC-class with envelope $F$. Let $\mathscr{Q}$ be the set of finite discrete probability measures $Q$ on $\mathbb{H}_+$ such that $\|F\|_{Q,2} > 0$. For $Q \in \mathscr{Q}$ and a Borel-measurable function $f : \mathbb{H}_+ \to \mathbb{R}$ satisfying $\|f\|_{Q,2} < \infty$, let $B_f(\varepsilon) = \{g \in \mathscr{V} : \|f - g\|_{Q,2} < \varepsilon\}$ denote the $L_2(Q)$-ball in $\mathscr{V}$ with center $f$ and radius $\varepsilon$. Let $N(\varepsilon, \mathscr{V}, L_2(Q))$ be the smallest number of balls $B_f(\varepsilon)$ needed to cover $\mathscr{V}$. Then $\mathscr{V}$ satisfies

$$(4.5) \qquad \int_0^\infty \sup_{Q \in \mathscr{Q}} \sqrt{\log N(\varepsilon \|F\|_{Q,2}, \mathscr{V}, L_2(Q))} \, d\varepsilon < \infty.$$

See Definition 2.1.5, (2.5.1) and Theorem 2.6.7 in [23].

The uniform entropy condition (4.5) is used in the next lemma to show that, for large $k \in \mathbb{N}$, the $\ell^\infty(\mathscr{V})$-norm of the random measure

$$(4.6) \qquad \hat{\mathcal{G}}_k^r = \frac{1}{\sqrt{k}} \sum_{i=1}^k (\delta_{(v_i^r, l_i^r r^{-1})} - \breve{\mathscr{v}}^r)$$

satisfies a tail bound that is uniform in $r \in \mathcal{R}$. In particular, (4.5) implies that for each $r \in \mathcal{R}$, the class $\mathscr{V}$ admits a $\breve{\mathscr{v}}^r$-*Brownian bridge*, that is, a tight, Borel-measurable version of the zero-mean Gaussian process $\{\mathcal{G}_*^r(f) : f \in \mathscr{V}\}$ with covariance function

$$\mathbf{E}^r[\mathcal{G}_*^r(f)\mathcal{G}_*^r(g)] = \langle fg, \breve{\mathscr{v}}^r \rangle - \langle f, \breve{\mathscr{v}}^r \rangle \langle g, \breve{\mathscr{v}}^r \rangle.$$

Moreover, the expected norms $\mathbf{E}^r[\|\langle f, \mathcal{G}_*^r \rangle\|_{\mathscr{V}}]$ are uniformly bounded in $r$ and for large $k$, the random measures $\{\hat{\mathcal{G}}_k^r : r \in \mathcal{R}\}$ are uniformly close in $\ell^\infty(\mathscr{V})$ to $\{\mathcal{G}_*^r : r \in \mathcal{R}\}$. This quickly leads to the desired tail bound; see Lemma 4.6 below. Finally, Lemma 4.7 below employs this tail bound to establish the functional Glivenko–Cantelli estimate for $\bar{\mathcal{L}}_{s,t}^{r,m}$.

The proofs in this section use outer probabilities and outer expectations to avoid cumbersome verifications of the measurability of certain suprema over uncountable sets. For each $r \in \mathcal{R}$ and each subset $S \subset \Omega^r$, define the *outer probability*

$$\dot{\mathbf{P}}^r(S) = \inf\{\mathbf{P}^r(B) : S \subset B, B \in \mathscr{F}^r\}$$

and the *outer expectation* of a map $X : \Omega^r \to \mathbb{R}$ by

$$\dot{\mathbf{E}}^r[X] = \inf\{\mathbf{E}^r[Y] : X \leq Y, Y \text{ is } \mathscr{F}^r\text{-measurable}\}.$$

Assume (A) and let $\mathscr{V}$ be a VC-class of Borel-measurable functions $f : \mathbb{H}_+ \to \mathbb{R}$ with Borel-measurable envelope $F$ such that $\mathscr{V}_\infty^2$ and $\mathscr{V}_\delta^r$ are Borel-measurable classes for all $r \in \mathcal{R}$ and $\delta > 0$ and such that

$$(4.7) \qquad \lim_{N \to \infty} \sup_{r \in \mathcal{R}} \langle F^2 1_{\{F > N\}}, \breve{\mathscr{v}}^r \rangle = 0.$$

LEMMA 4.6. *For all $q \geq 1$ and $x > 2$, there exist $M < \infty$ and $k_0 \in \mathbb{N}$ such that $k > k_0$ implies*

$$\sup_{r \in \mathcal{R}} \dot{\mathbf{P}}^r(\|\langle f, \hat{\mathcal{G}}_k^r \rangle\|_{\mathscr{V}} > x) \leq \frac{M}{x^q}.$$

*The constant $M$ does not depend on $x$.*

PROOF. Fix $q \geq 1$ and $x > 2$. For $\mathcal{G} \in \ell^\infty(\mathscr{V})$, define

$$(4.8) \qquad h_x(\mathcal{G}) = (\|\mathcal{G}(f)\|_{\mathscr{V}} - x + 1)^+ \wedge 1.$$

Then

$$(4.9) \qquad \sup_{r \in \mathcal{R}} \dot{\mathbf{P}}^r(\|\langle f, \hat{\mathcal{G}}_k^r \rangle\|_{\mathscr{V}} > x) \leq \sup_{r \in \mathcal{R}} \dot{\mathbf{E}}^r[h_x(\hat{\mathcal{G}}_k^r)].$$

The function $h_x$ is an element of $\mathbf{BL_1}(\mathscr{V})$, the set of all $h : \ell^\infty(\mathscr{V}) \to \mathbb{R}$ satisfying $\|h\|_\infty \leq 1$ and $|h(\mathcal{G}_1) - h(\mathcal{G}_2)| \leq \|\mathcal{G}_1 - \mathcal{G}_2\|_{\mathscr{V}}$ for all $\mathcal{G}_1, \mathcal{G}_2 \in \ell^\infty(\mathscr{V})$. Since $\mathscr{V}$ is a VC-class, it satisfies the uniform entropy condition (4.5). Moreover, the envelope $F$ satisfies (4.7) and $\mathscr{V}_\infty^2, \mathscr{V}_\delta^r$ satisfy the indicated measurability conditions. Thus, $\mathscr{V}$ satisfies the assumptions of Theorem 2.8.3 in [23] which asserts that $\mathscr{V}$ is *pre-Gaussian* and *Donsker*, uniformly in $r \in \mathcal{R}$ (see Section 2.8.2 in [23]). In particular, $\mathscr{V}$ satisfies

$$(4.10) \qquad \sup_{r \in \mathcal{R}} \mathbf{E}^r[\|\langle f, \mathcal{G}_*^r \rangle\|_{\mathscr{V}}] < \infty$$

and since $h_x \in \mathbf{BL_1}$,

$$(4.11) \qquad \lim_{k \to \infty} \sup_{r \in \mathcal{R}} |\dot{\mathbf{E}}^r[h_x(\hat{\mathcal{G}}_k^r)] - \mathbf{E}^r[h_x(\mathcal{G}_*^r)]| = 0.$$

By (4.9) and (4.11), there exists $k_0 \in \mathbb{N}$ such that $k > k_0$ implies

$$\sup_{r \in \mathcal{R}} \dot{\mathbf{P}}^r(\|\langle f, \hat{\mathcal{G}}_k^r \rangle\|_{\mathscr{V}} > x) \leq \sup_{r \in \mathcal{R}} \mathbf{E}^r[h_x(\mathcal{G}_*)] + x^{-q}.$$

Apply (4.8) and Markov's inequality to the right-hand side to obtain

$$\sup_{r \in \mathcal{R}} \dot{\mathbf{P}}^r(\|\langle f, \hat{\mathcal{G}}_k^r \rangle\|_{\mathscr{V}} > x) \leq \sup_{r \in \mathcal{R}} \mathbf{P}^r(\|\langle f, \mathcal{G}_*^r \rangle\|_{\mathscr{V}} > x - 1) + x^{-q}$$

$$\leq x^{-q}\left(2^q \sup_{r \in \mathcal{R}} \mathbf{E}^r[\|\langle f, \mathcal{G}_*^r \rangle\|_{\mathscr{V}}^q] + 1\right).$$

Let $M$ be the last term in parentheses, which does not depend on $x$. For each $r \in \mathcal{R}$, the Brownian bridge $\mathcal{G}_*^r$ is separable and Gaussian with $\|\langle f, \mathcal{G}_*^r \rangle\|_{\mathscr{V}}$ finite almost surely. Thus, there exists a constant $C$ such that for all $r \in \mathcal{R}$,

$$\mathbf{E}^r[\|\langle f, \mathcal{G}_*^r \rangle\|_{\mathscr{V}}^q] \leq C\mathbf{E}^r[\|\langle f, \mathcal{G}_*^r \rangle\|_{\mathscr{V}}]^q$$

(see Proposition A.2.4 in [23]). Conclude from (4.10) that $M < \infty$. $\quad\square$

The next lemma establishes the functional Glivenko–Cantelli estimate for $\bar{\mathcal{L}}_{s,t}^{r,m}$. Assume (A) and let $T, L > 1$ and $\varepsilon \in (0,1)$. Let $\mathscr{V}$ be a VC-class of Borel-measurable functions $f \colon \mathbb{H}_+ \to \mathbb{R}$ such that $\mathscr{V}_\infty^2$ and $\mathscr{V}_\delta^r$ are Borel-measurable classes for all $r \in \mathcal{R}$ and $\delta > 0$ and such that (4.7) holds.

LEMMA 4.7. *There exist events*

$$\Omega_1^r \subset \left\{ \sup_{f \in \mathscr{V}} \sup_{\substack{m \leq \lfloor rT \rfloor, \\ [s,t] \subset [0,L]}} |\langle f, \bar{\mathcal{L}}_{s,t}^{r,m}\rangle - \alpha^r(t-s)\langle f, \breve{\vartheta}^r\rangle| \leq \varepsilon \right\}, \qquad r \in \mathcal{R},$$

*such that* $\lim_{r \to \infty} \mathbf{P}^r(\Omega_1^r) = 1$.

PROOF. For each $r \in \mathcal{R}$, let

$$(4.12) \qquad \Upsilon^r = \left\{ \sup_{\substack{m \leq \lfloor rT \rfloor, \\ [s,t] \subset [0,L]}} \|\langle f, \bar{\mathcal{L}}_{s,t}^{r,m}\rangle - \alpha^r(t-s)\langle f, \breve{\vartheta}^r\rangle\|_{\mathscr{V}} > \varepsilon \right\}.$$

It suffices to show that

$$(4.13) \qquad \qquad \lim_{r \to \infty} \mathring{\mathbf{P}}^r(\Upsilon^r) = 0.$$

Since

$$\langle f, \bar{\mathcal{L}}_{s,t}^{r,m}\rangle - \alpha^r(t-s)\langle f, \breve{\vartheta}^r\rangle = (\langle f, \bar{\mathcal{L}}^r(m+t)\rangle - \alpha^r(m+t)\langle f, \breve{\vartheta}^r\rangle)$$
$$- (\langle f, \bar{\mathcal{L}}^r(m+s)\rangle - \alpha^r(m+s)\langle f, \breve{\vartheta}^r\rangle),$$

the supremum in (4.12) is bounded above by

$$\sup_{\substack{[s,t] \subset [0, rT+L], \\ |t-s| \leq L}} \|(\langle f, \bar{\mathcal{L}}^r(t)\rangle - \alpha^r t\langle f, \breve{\vartheta}^r\rangle) - (\langle f, \bar{\mathcal{L}}^r(s)\rangle - \alpha^r s\langle f, \breve{\vartheta}^r\rangle)\|_{\mathscr{V}}.$$

Rewrite this as

$$\sup_{\substack{[s,t] \subset [0, T+L/r], \\ |t-s| \leq L/r}} \|(\langle f, \bar{\mathcal{L}}^r(rt)\rangle - \alpha^r rt\langle f, \breve{\vartheta}^r\rangle) - (\langle f, \bar{\mathcal{L}}^r(rs)\rangle - \alpha^r rs\langle f, \breve{\vartheta}^r\rangle)\|_{\mathscr{V}}.$$

For $u \geq 0$, let $\bar{\bar{E}}^r(u) = r^{-2}E^r(r^2u)$ and $\hat{E}^r(u) = \bar{E}^r(ru) - \alpha^r ru$. Then

$$\langle f, \bar{\mathcal{L}}^r(ru)\rangle - \alpha^r ru\langle f, \breve{\vartheta}^r\rangle = \frac{1}{r}\sum_{i=1}^{r\bar{E}^r(ru)} (f(v_i^r, l_i^r r^{-1}) - \langle f, \breve{\vartheta}^r\rangle)$$
$$+ (\bar{E}^r(ru) - \alpha^r ru)\langle f, \breve{\vartheta}^r\rangle$$
$$= \bar{\bar{E}}^r(u)^{1/2}\langle f, \hat{\mathcal{G}}_{r^2\bar{E}^r(u)}^r\rangle + \hat{E}^r(u)\langle f, \breve{\vartheta}^r\rangle.$$

Introducing the notation $X_k^r(f) = (k/r^2)^{1/2}\langle f, \hat{\mathcal{G}}_k^r\rangle$ for $k \in \mathbb{N}$ and $f \in \mathcal{V}$, the supremum in (4.12) is bounded above by

$$\sup_{\substack{[s,t] \subset [0, T+L/r], \\ |t-s| \le L/r}} \|X_{r^2\bar{E}^r(t)}^r(f) - X_{r^2\bar{E}^r(s)}^r(f) + (\hat{E}^r(t) - \hat{E}^r(s))\langle f, \breve{\vartheta}^r\rangle\|_{\mathcal{V}}.$$

Let $\delta > 0$. Deduce from the previous remarks and the inequalities $\|\langle f, \breve{\vartheta}^r\rangle\|_{\mathcal{V}} \le \langle F, \breve{\vartheta}^r\rangle$ and $\limsup_{r\to\infty} L/r \le T \wedge \delta$ that

$$\lim_{r\to\infty} \dot{\mathbf{P}}(\Upsilon^r) \le \limsup_{r\to\infty} \dot{\mathbf{P}}^r\left(\sup_{\substack{s,t \le 2T, \\ |t-s| \le \delta}} |\hat{E}^r(t) - \hat{E}^r(s)|\langle f, \breve{\vartheta}^r\rangle > \frac{\varepsilon}{2}\right)$$

(4.14)

$$+ \limsup_{r\to\infty} \dot{\mathbf{P}}^r\left(\sup_{\substack{s,t \le 2T, \\ |t-s| \le \delta}} \|X_{r^2\bar{E}^r(t)}^r(f) - X_{r^2\bar{E}^r(s)}^r(f)\|_{\mathcal{V}} > \frac{\varepsilon}{2}\right).$$

By (2.22)–(2.24) and the functional central limit theorem for renewal processes, the sequence $\{\hat{E}^r(\cdot)\}$ converges in distribution to a Brownian motion which is almost surely uniformly continuous on compact time intervals [2]. Since $\sup_{r \in \mathcal{R}} \langle F, \breve{\vartheta}^r\rangle < \infty$ by (4.7), we can conclude that the first right-hand term in (4.14) converges to zero as $\delta \to 0$. To show the same for the second term, the proof of Theorem 2.12.1 in [23] has been adapted.

For $a, b \in [0, \infty)$, define $\mathbb{N}_a^b = \mathbb{N} \cap [a, b]$. It suffices to show that

$$(4.15) \quad \limsup_{r\to\infty} \dot{\mathbf{P}}^r\left(\max_{j \in \mathbb{N}_1^{\lceil 2T/\delta\rceil}} \sup_{s,t \in [(j-1)\delta, j\delta]} \|X_{r^2\bar{E}^r(t)}^r(f) - X_{r^2\bar{E}^r(s)}^r(f)\|_{\mathcal{V}} > \frac{\varepsilon}{4}\right)$$

converges to zero as $\delta \to 0$. For each $r \in \mathcal{R}$, bound the outer probability of the maximum by the sum of outer probabilities and discretize the time index in the supremum. Then (4.15) is bounded above by

$$(4.16) \quad \limsup_{r\to\infty} \sum_{j=1}^{\lceil 2T/\delta\rceil} \dot{\mathbf{P}}^r\left(\max_{k,l \in \mathbb{N}_{r^2\bar{E}^r((j-1)\delta)}^{r^2\bar{E}^r(j\delta)}} \|X_l^r(f) - X_k^r(f)\|_{\mathcal{V}} > \frac{\varepsilon}{4}\right).$$

Let $\Omega_0^r = \{\sup_{u\in[0,2T+\delta]} |\bar{\bar{E}}^r(u) - \alpha^r u| > \alpha\delta\}$. Then on the complement of $\Omega_0^r$, for each $j = 1, \dots, \lceil 2T/\delta\rceil$,

$$(\alpha^r j\delta - \alpha\delta)^+ \le \bar{\bar{E}}^r(j\delta) \le \alpha^r j\delta + \alpha\delta.$$

Deduce that (4.16) is bounded above by

$$\limsup_{r\to\infty} \sum_{j=1}^{\lceil 2T/\delta\rceil} \dot{\mathbf{P}}^r\left(\max_{k,l \in \mathbb{N}_{r^2(\alpha^r(j-1)\delta-\alpha\delta)^+}^{r^2(\alpha^r j\delta+\alpha\delta)}} \left\|\frac{\sqrt{l}}{r}\langle f, \hat{\mathcal{G}}_l^r\rangle - \frac{\sqrt{k}}{r}\langle f, \hat{\mathcal{G}}_k^r\rangle\right\|_{\mathcal{V}} > \frac{\varepsilon}{4}\right)$$

$$+ \limsup_{r\to\infty} \left\lceil\frac{2T}{\delta}\right\rceil \dot{\mathbf{P}}^r(\Omega_0^r).$$

The second term equals zero by (2.22)–(2.24) and the functional weak law of large numbers for renewal processes. Since $\hat{\mathcal{G}}_k^r$ has stationary increments as a function of $k$ and $\limsup_{r\to\infty}\alpha^r < 2\alpha$, the previous term is bounded above by

$$(4.17) \qquad \limsup_{r\to\infty}\left\lceil\frac{2T}{\delta}\right\rceil\dot{\mathbf{P}}^r\left(\max_{k\le\lfloor r^2 4\alpha\delta\rfloor}\left\|\frac{\sqrt{k}}{r}\langle f,\hat{\mathcal{G}}_k^r\rangle\right\|_{\mathscr{V}} > \frac{\varepsilon}{4}\right).$$

By Ottaviani's inequality (see Proposition A.1.1 in [23]) and the stationary increments of $\hat{\mathcal{G}}^r$, (4.17) is bounded above by

$$(4.18) \qquad \limsup_{r\to\infty}\frac{\lceil 2T/\delta\rceil\dot{\mathbf{P}}^r(\|\langle f,\hat{\mathcal{G}}_{\lfloor r^2 4\alpha\delta\rfloor}^r\rangle\|_{\mathscr{V}} > \varepsilon/(16\sqrt{\alpha\delta}))}{1 - \max_{k\le\lfloor r^2 4\alpha\delta\rfloor}\dot{\mathbf{P}}^r(\|\langle f,\hat{\mathcal{G}}_k^r\rangle\|_{\mathscr{V}} > \varepsilon r/(8\sqrt{k}))}.$$

Assume that $\delta$ is sufficiently small that $\varepsilon(16\sqrt{\alpha\delta})^{-1} > 2$. By Lemma 4.6, there exist $M$ and $k_0\in\mathbb{N}$ such that $k > k_0$ implies

$$(4.19) \qquad \sup_{r\in\mathcal{R}}\dot{\mathbf{P}}^r\left(\|\langle f,\hat{\mathcal{G}}_k^r\rangle\|_{\mathscr{V}} > \frac{\varepsilon}{16\sqrt{\alpha\delta}}\right) \le \left(\frac{16\sqrt{\alpha\delta}}{\varepsilon}\right)^3 M.$$

Since $\lfloor r^2 4\alpha\delta\rfloor\to\infty$ as $r\to\infty$, the limit superior of the numerator in (4.18) is bounded above by $\lceil 2T/\delta\rceil(16\sqrt{\alpha\delta}/\varepsilon)^3 M$, which can be made arbitrarily small by choosing $\delta$ sufficiently close to 0. It remains to show that the limit inferior, as $r\to\infty$, of the denominator in (4.18) is bounded away from zero as $\delta\to 0$. Observe that $\varepsilon r(8\sqrt{k})^{-1}\ge\varepsilon(16\sqrt{\alpha\delta})^{-1}$ for $k\le\lfloor r^2 4\alpha\delta\rfloor$. So, by (4.19), the terms in the maximum indexed by $k > k_0$ are bounded above (uniformly in $r$) by $(16\sqrt{\alpha\delta}/\varepsilon)^3 M$, which can be made arbitrarily small by choosing $\delta$ sufficiently close to 0. For $k\le k_0$,

$$\dot{\mathbf{P}}^r\left(\|\langle f,\hat{\mathcal{G}}_k^r\rangle\|_{\mathscr{V}} > \frac{\varepsilon r}{8\sqrt{k}}\right) \le \dot{\mathbf{P}}^r\left(\|\langle f,\hat{\mathcal{G}}_k^r\rangle\|_{\mathscr{V}} > \frac{\varepsilon r}{8\sqrt{k_0}}\right),$$

which converges to zero as $r\to\infty$. Conclude that (4.18) converges to zero as $\delta\to 0$.  □

It remains to define a class of subsets of $\mathbb{H}_+$ that is sufficiently rich to characterize elements of $\mathbf{M}$, yet sufficiently small to be a VC-class.

DEFINITION 4.8.   Let $\overline{(a,b)}$ denote the closure in $\mathbb{R}$ of an interval $(a,b)\subset\mathbb{R}$ and define

$$\mathscr{A} = \{[x,\infty)\times\overline{(y,\infty)}\colon x\in[0,\infty), y\in[-\infty,\infty)\}.$$

Note that since the above definition allows $y = -\infty$, the collection $\mathscr{A}$ includes not only subsets of $\mathbb{H}_+$ that are translations of the right upper quadrant, but also all right half-spaces.

LEMMA 4.9.    *The family $\mathscr{V} = \{1_A : A \in \mathscr{A}\}$ is a VC-class of Borel-measurable functions with Borel-measurable envelope $F \equiv 1$. Moreover, $\mathscr{V}_\infty^2$ and $\mathscr{V}_\delta^r$ are Borel-measurable classes for all $r \in \mathcal{R}$ and $\delta > 0$ and the envelope $F$ satisfies* (4.7).

PROOF.    Since $\mathscr{V}$ contains only indicator functions, it suffices to show that $\mathscr{A}$ is a VC-class. Let $\{a, b, c\} \subset \mathbb{H}_+$ be a three-point subset. It is impossible for $\mathscr{A}$ to pick out all three two-point subsets of $\{a, b, c\}$. Thus, $\mathscr{A}$ shatters no three-point subset of $\mathbb{H}_+$ and so has VC-index at most 3. The remaining assertions are evident.    $\square$

4.3. *Dynamic equation.*    This section introduces the dynamic equation satisfied by the measure valued state descriptor. For a subset $B \subset \mathbb{H}_+$ and $w \in \mathbb{H}_+$, define $1_B^+(w) = \langle 1_B, \delta_w^+ \rangle$ (recall the definition of $\delta_w^+$ from Section 2.1). For each $r \in \mathcal{R}$, the state descriptor of the $r$th model satisfies the following equation almost surely: for each Borel set $B \subset \mathbb{H}_+$ and all $t, h \geq 0$,

$$\mathcal{Z}^r(t+h)(B) = \mathcal{Z}^r(t)(B + (S_{t,t+h}^r, h))$$
$$+ \sum_{i=E^r(t)+1}^{E^r(t+h)} 1_B^+(v_i^r(t+h), l_i^r(t+h)).$$

This follows from (2.3)–(2.8) after some simplification. Applying the fluid scaling (2.10)–(2.15) and (2.17), the above equation becomes

(4.20)
$$\bar{\mathcal{Z}}^r(t+h)(B) = \bar{\mathcal{Z}}^r(t)(B + (\bar{S}_{t,t+h}^r, h))$$
$$+ \frac{1}{r} \sum_{i=r\bar{E}^r(t)+1}^{r\bar{E}^r(t+h)} 1_B^+(\bar{v}_i^r(t+h), \bar{l}_i^r(t+h)).$$

After rescaling and shifting, the equation takes the following form. For each $r \in \mathcal{R}$, $m \leq \lfloor rT \rfloor$ and $t, h \geq 0$, almost surely for each Borel set $B \subset \mathbb{H}_+$,

(4.21)
$$\bar{\mathcal{Z}}^{r,m}(t+h)(B) = \bar{\mathcal{Z}}^{r,m}(t)(B + (\bar{S}_{t,t+h}^{r,m}, h))$$
$$+ \frac{1}{r} \sum_{i=r\bar{E}^{r,m}(t)+1}^{r\bar{E}^{r,m}(t+h)} 1_B^+(\bar{v}_i^{r,m}(t+h), \bar{l}_i^{r,m}(t+h)).$$

Equations (4.20) and (4.21) are called the *dynamic equations* for $\bar{\mathcal{Z}}^r(\cdot)$ and $\bar{\mathcal{Z}}^{r,m}(\cdot)$.

Subsequent proofs use estimates obtained from (4.21). Two estimates result from bounding the summands by 1 and optionally bounding the first term on the right-hand side by its total mass; for each $B \subset \mathbb{H}_+$ and $t, h \geq 0$,

$$(4.22) \quad \begin{aligned} \bar{\mathcal{Z}}^{r,m}(t+h)(B) &\leq \bar{\mathcal{Z}}^{r,m}(t)(B + (\bar{S}^{r,m}_{t,t+h}, h)) + \bar{\mathcal{L}}^{r,m}_{t,t+h}(\mathbb{H}_+) \\ &\leq \bar{\mathcal{Z}}^{r,m}(t)(\mathbb{H}_+) + \bar{\mathcal{L}}^{r,m}_{t,t+h}(\mathbb{H}_+). \end{aligned}$$

Two more estimates follow from (4.21) by simply ignoring any arrivals; for each $B \subset \mathbb{H}_+$ and $t, h \geq 0$,

$$(4.23) \quad \bar{\mathcal{Z}}^{r,m}(t)(B + (\bar{S}^{r,m}_{t,t+h}, h)) \leq \bar{\mathcal{Z}}^{r,m}(t+h)(B) \leq \bar{\mathcal{Z}}^{r,m}(t+h)(\mathbb{H}_+).$$

4.4. *Compact containment.* The following lemma establishes the compact containment of the fluid scaled state descriptor on $[0, \lfloor rT \rfloor + L]$. Assume (A) and let $T > 1$ and $\eta \in (0, 1)$.

LEMMA 4.10.   *There exists a compact set* $\mathbf{K} \subset \mathbf{M}$ *such that for all* $L > 1$,

$$\liminf_{r \to \infty} \mathbf{P}^r(\bar{\mathcal{Z}}^r(t) \in \mathbf{K}) \geq 1 - \eta \quad \text{for all } t \in [0, \lfloor rT \rfloor + L]) \geq 1 - \eta.$$

PROOF.   Since $\pi_+ \colon (x, y) \mapsto |y|$ is continuous, (2.25) and the Skorohod representation theorem imply the existence of $\mathbb{R}_+$-valued random variables $X^r \sim \breve{\vartheta}^r \circ \pi_+^{-1}$ and $X \sim \vartheta \circ \pi_+^{-1}$ such that $X^r \to X$ almost surely. Consequently, there exists an $\mathbb{R}_+$-valued random variable $Y$ such that

$$(4.24) \qquad\qquad Y = \sup_r X^r \qquad \text{a.s.}$$

Let $\mu$ be the law of $Y$ on $\mathbb{R}_+$. Since $L^2(\mu)$ contains unbounded functions, there exists a continuous increasing unbounded $\psi \colon \mathbb{R}_+ \to \mathbb{R}_+$ such that $\langle \psi^2, \mu \rangle < \infty$. This implies that

$$\langle (\psi \circ \pi_+)^2, \vartheta \rangle = \mathbf{E}[\psi(X)^2] \leq \mathbf{E}[\psi(Y)^2] < \infty.$$

Thus, the one-element set $\mathscr{V} = \{\psi \circ \pi_+\}$ is a VC-class with envelope $F = \psi \circ \pi_+$ satisfying (4.7).

Since $\mathbf{M}$ is a Polish space, (2.32) and Prohorov's theorem imply that the sequence $\{\bar{\mathcal{Z}}^r(0) \colon r \in \mathcal{R}\}$ is tight: there exists a compact $\mathbf{K}_0 \subset \mathbf{M}$ such that

$$(4.25) \qquad\qquad \liminf_{r \to \infty} \mathbf{P}^r(\bar{\mathcal{Z}}^r(0) \in \mathbf{K}_0) \geq 1 - \frac{\eta}{3}.$$

By (4.25), Corollary 4.3 and Lemma 4.7, there exist $M > 1$ and events

$$\Omega^r_1 \subset \left\{ \sup_{m \leq \lfloor 2rT \rfloor} |\langle \psi \circ \pi_+, \bar{\mathcal{L}}^{r,m}_{0,1} \rangle - \alpha^r \langle \psi \circ \pi_+, \breve{\vartheta}^r \rangle | \leq 1 \right\}, \qquad r \in \mathcal{R},$$

such that for all $L > 1$, the events

$$\Omega^r_2 = \{ \bar{\mathcal{Z}}^r(0) \in \mathbf{K}_0 \},$$

$$\Omega^r_3 = \left\{ \sup_{t \in [0, \lfloor rT \rfloor + L]} \langle 1, \bar{\mathcal{Z}}^r(t) \rangle \vee \langle \chi, \bar{\mathcal{Z}}^r(t) \rangle \leq M \right\},$$

$$\Omega^r_0 = \Omega^r_1 \cap \Omega^r_2 \cap \Omega^r_3$$

satisfy

(4.26)                        $\liminf_{r\to\infty} \mathbf{P}^r(\Omega_0^r) \geq 1 - \eta.$

Assume henceforth in the proof that all random objects are evaluated at some outcome $\omega \in \Omega_0^r$.

Since $\mathbf{K}_0$ is compact, there exists a sequence of real numbers $\{x_k : k \in \mathbb{N}\}$ such that for all $k$,

(4.27)                        $\sup_{\xi \in \mathbf{K}_0} \xi(([0, x_k] \times [-x_k, x_k])^c) \leq \frac{1}{k}.$

For each $k \in \mathbb{N}$, choose $y_k \geq x_k$ such that $\psi(y_k) \geq k^2$ and define

$$I_k = [0, k],$$

$$J_k = [-(y_k + k), y_k],$$

$$R_k = I_k \times J_k.$$

Consider a sequence $a_k \to 0$ as $k \to \infty$. The set

(4.28)                $\mathbf{K} = \{\xi \in \mathbf{M} : \langle 1, \xi \rangle \leq M \text{ and } \xi(R_k^c) \leq a_k \text{ for all } k\}$

is precompact ([12], Theorem A 7.5). Thus, by (4.26) and the definition of $\Omega_3^r$, it suffices to identify a null sequence $\{a_k\}$ such that for each $L > 1$ and each $r > L$,

$$\sup_{t \in [0, \lfloor rT \rfloor + L]} \bar{\mathcal{Z}}^r(t)(R_k^c) \leq a_k \qquad \text{for all } k.$$

Fix $L > 1$, $t \in [0, \lfloor rT \rfloor + L]$ and $k \in \mathbb{N}$. Assume that $r > L$. Observe that

(4.29)                $\bar{\mathcal{Z}}^r(t)(R_k^c) = \bar{\mathcal{Z}}^r(t)(I_k^c \times \mathbb{R}) + \bar{\mathcal{Z}}^r(t)(I_k \times J_k^c).$

Recall that $\chi(x, y) = x$. By Markov's inequality and the definition of $\Omega_3^r$,

(4.30)        $\bar{\mathcal{Z}}^r(t)(I_k^c \times \mathbb{R}) = \bar{\mathcal{Z}}^r(t)((k, \infty) \times \mathbb{R}) \leq \frac{1}{k}\langle \chi, \bar{\mathcal{Z}}^r(t) \rangle \leq \frac{M}{k}.$

Let $\mathcal{N} = \{s \in [(t-k)^+, t] : \langle 1, \bar{\mathcal{Z}}^r(s) \rangle = 0\}$; if $\mathcal{N} \neq \varnothing$, fix some $\tau \in \mathcal{N}$ and if $\mathcal{N} = \varnothing$, let $\tau = (t-k)^+$. By (4.20),

(4.31)          $\bar{\mathcal{Z}}^r(t)(I_k \times J_k^c) = \bar{\mathcal{Z}}^r(\tau)((I_k \times J_k^c) + (\bar{S}_{\tau,t}^r, t - \tau))$

$$+ \frac{1}{r}\sum_{i=r\bar{E}^r(\tau)+1}^{r\bar{E}^r(t)} 1^+_{I_k \times J_k^c}(\bar{v}_i^r(t), \bar{l}_i^r(t)).$$

If $\mathcal{N} \neq \varnothing$, then $\bar{\mathcal{Z}}^r(\tau) = \mathbf{0}$ and

(4.32)              $\bar{\mathcal{Z}}^r(\tau)((I_k \times J_k^c) + (\bar{S}_{\tau,t}^r, t - \tau)) = 0.$

If $\mathcal{N} = \varnothing$, then either $\tau = 0$ or $\tau = t - k$. Suppose $\tau = 0$. Then by the definitions of $\Omega_2^r$, $y_k$ and $J_k$, by (4.27) and by the inequality $t - \tau \leq k$,

$$(4.33) \quad \bar{\mathcal{Z}}^r(\tau)((I_k \times J_k^c) + (\bar{S}_{\tau,t}^r, t - \tau)) \leq \bar{\mathcal{Z}}^r(0)(\mathbb{R}_+ \times [-y_k, y_k]^c) \leq \frac{1}{k}.$$

Suppose $\tau = t - k$. Then, by definition of $\Omega_3^r$ and the fact that $\mathcal{N} = \varnothing$, we have

$$\bar{S}_{\tau,t}^r = \int_\tau^t \langle 1, \bar{\mathcal{Z}}^r(s) \rangle^{-1} \, ds \geq \frac{t - \tau}{M} = \frac{k}{M}.$$

In this case, Markov's inequality and the definition of $\Omega_3^r$ yield

$$(4.34) \quad \begin{aligned} \bar{\mathcal{Z}}^r(\tau)((I_k \times J_k^c) + (\bar{S}_{\tau,t}^r, t - \tau)) &\leq \bar{\mathcal{Z}}^r(\tau)([k/M, \infty) \times \mathbb{R}) \\ &\leq \frac{M}{k} \langle \chi, \bar{\mathcal{Z}}^r(\tau) \rangle \\ &\leq \frac{M^2}{k}. \end{aligned}$$

We can then deduce from (4.31)–(4.34) that

$$(4.35) \quad \bar{\mathcal{Z}}^r(t)(I_k \times J_k^c) \leq \frac{M^2}{k} + \frac{1}{r} \sum_{i=r\bar{E}^r(\tau)+1}^{r\bar{E}^r(t)} 1_{I_k \times J_k^c}^+ (\bar{v}_i^r(t), \bar{l}_i^r(t)).$$

For each $i$ appearing in the sum, (2.5) and (2.15) imply that

$$\bar{l}_i^r(t) = l_i^r r^{-1} - t + U_i^r r^{-1} \geq l_i^r r^{-1} - t + \tau \geq l_i^r r^{-1} - k.$$

Thus, $\bar{l}_i^r(t) \in [-(y_k + k), y_k]^c$ implies that $l_i^r r^{-1} \in [-y_k, y_k]^c$, yielding

$$1_{I_k \times J_k^c}^+ (\bar{v}_i^r(t), \bar{l}_i^r(t)) \leq 1_{\mathbb{R}_+ \times [-y_k, y_k]^c}(v_i^r, l_i^r r^{-1}).$$

By definition of $\psi$ and $y_k$,

$$1_{\mathbb{R}_+ \times [-y_k, y_k]^c}(v_i^r, l_i^r r^{-1}) \leq \frac{1}{\psi(y_k)} \psi(|l_i^r r^{-1}|) \leq \frac{1}{k^2} \psi(|l_i^r r^{-1}|).$$

Deduce from (4.35) that

$$\bar{\mathcal{Z}}^r(t)(I_k \times J_k^c) \leq \frac{M^2}{k} + \frac{1}{r} \sum_{i=r\bar{E}^r(\tau)+1}^{r\bar{E}^r(t)} \frac{1}{k^2} \psi(|l_i^r r^{-1}|).$$

Bounding and rewriting the sum, we obtain

$$\bar{\mathcal{Z}}^r(t)(I_k \times J_k^c) \leq \frac{M^2}{k} + \frac{1}{k^2} \sum_{m=\lfloor \tau \rfloor}^{\lfloor t \rfloor} \frac{1}{r} \sum_{i=r\bar{E}^{r,m}(0)+1}^{r\bar{E}^{r,m}(1)} \psi(|l_i^r r^{-1}|).$$

By assumption, $t \leq rT + L = r(T + Lr^{-1}) \leq 2rT$. Rewriting the sum again and using the definition of $\Omega_1^r$, we get

$$
\begin{aligned}
\bar{\mathcal{Z}}^r(t)(I_k \times J_k^c) &\leq \frac{M^2}{k} + \frac{1}{k^2} \sum_{m=\lfloor \tau \rfloor}^{\lfloor t \rfloor} \langle \psi \circ \pi_+, \bar{\mathcal{L}}_{0,1}^{r,m} \rangle \\
(4.36) \\
&\leq \frac{M^2}{k} + \frac{k+1}{k^2}(\alpha^r \langle \psi \circ \pi_+, \breve{\vartheta}^r \rangle + 1).
\end{aligned}
$$

By (2.22), (4.24) and the definition of $\psi$, there exists $K < \infty$ such that

$$
(4.37) \quad \sup_{r \in \mathcal{R}} \alpha^r \langle \psi \circ \pi_+, \breve{\vartheta}^r \rangle = \sup_{r \in \mathcal{R}} \alpha^r \mathbf{E}[\psi(X^r)] \leq \sup_{r \in \mathcal{R}} \alpha^r \mathbf{E}[\psi(Y)] \leq K.
$$

Combining (4.29), (4.30), (4.36) and (4.37) yields

$$
(4.38) \quad \bar{\mathcal{Z}}^r(t)(R_k^c) \leq \frac{1}{k}(M + M^2 + 2(K+1)).
$$

Let $a_k$ equal the right-hand side of (4.38). Since $a_k$ does not depend on $L$ and since (4.38) holds on $\Omega_0^r$ for all $L > 1$, $r > L$ and $t \in [0, \lfloor rT \rfloor + L]$, the proof is complete.  $\square$

4.5. *Asymptotic regularity.* This section establishes that, as $r \to \infty$ the fluid scaled state descriptors $\bar{\mathcal{Z}}^r(\cdot)$ assign arbitrarily small mass to the boundaries of sets $A \in \mathscr{A}$. This is phrased in terms of $\kappa$-enlargements of the boundaries of these sets. For $B \subset \mathbb{H}_+$ and $\kappa > 0$, let $\partial_B$ denote the boundary of $B$ and recall that $\partial_B^\kappa = \{w \in \mathbb{H}_+ : \inf_{z \in \partial_B} \|w - z\| < \kappa\}$ is the $\kappa$-enlargement in $\mathbb{H}_+$ of its boundary.

The result is established first for the initial condition $\bar{\mathcal{Z}}^r(0)$. Assume (A) and let $\varepsilon, \eta \in (0, 1)$.

LEMMA 4.11.  *There exists $\kappa > 0$ such that*

$$
(4.39) \quad \liminf_{r \to \infty} \mathbf{P}^r \left( \sup_{A \in \mathscr{A}} \bar{\mathcal{Z}}^r(0)(\partial_A^\kappa) \leq \varepsilon \right) \geq 1 - \eta.
$$

PROOF.  Note that the event in (4.39) is $\mathbf{P}^r$-measurable for each $r$; since the random Borel measure $\bar{\mathcal{Z}}^r(0)$ is continuous from below, the event can be rewritten using the supremum over a countable family $\{\partial_{A_n}^q\}$, where $q$ is rational and $\{A_n\} \subset \mathscr{A}$ have boundary points with rational coordinates.

By the first component of (2.32) and the Skorohod representation theorem, there exist random measures $\Lambda^r \sim \bar{\mathcal{Z}}^r(0)$ and $\Lambda \sim \Theta$, defined on a common probability space $(\Omega, \mathscr{F}, \mathbf{P})$, such that $\Lambda^r \xrightarrow{\mathbf{w}} \Lambda$ almost surely. It suffices to show (4.39) for $\Lambda^r$ in place of $\bar{\mathcal{Z}}^r(0)$. Since $\mathbf{M}$ is Polish, the sequence $\{\Lambda^r\}$ is strongly tight: there exists a compact $\mathbf{K} \subset \mathbf{M}$ such that $\Omega_0 = \bigcap_{r \in \mathcal{R}} \{\Lambda^r \in \mathbf{K}\}$ satisfies

$$
(4.40) \quad \mathbf{P}(\Omega_0) \geq 1 - \eta
$$

(see [17], Corollary 2). Let $\mathbf{K}_\vartheta = \mathbf{K} \cap \mathbf{M}_\vartheta$. By (2.30) and the fact that $\mathbf{K}$ is compact,

$$(4.41) \qquad\qquad \Lambda \in \mathbf{K}_\vartheta \qquad \text{a.s. on } \Omega_0.$$

Choose a compact $C \subset \mathbb{H}_+$ such that

$$\sup_{\xi \in \mathbf{K}} \xi(C^c) \leq \frac{\varepsilon}{2}$$

(see [12], Theorem A 7.5). For each $\kappa > 0$, let $\mathrm{I}^\kappa$ be the collection of all vertical strips $[i\kappa, (i+1)\kappa] \times \mathbb{R}$, $i = 0, 1, \ldots$, and all horizontal strips $\mathbb{R}_+ \times [j\kappa, (j+1)\kappa]$, $j \in \mathbb{Z}$. Define

$$\mathrm{I}_C^\kappa = \{I \in \mathrm{I}^\kappa : I \cap C \neq \varnothing\},$$

which is finite because $C$ is compact. Note that for each $A \in \mathscr{A}$, the set $\partial_A^\kappa \cap C$ is contained in the union of at most six strips in $\mathrm{I}_C^\kappa$. For all $\kappa > 0$ and $\xi \in \mathbf{K}$,

$$\sup_{A \in \mathscr{A}} \xi(\partial_A^\kappa) \leq \sup_{A \in \mathscr{A}} \xi(\partial_A^\kappa \cap C) + \frac{\varepsilon}{2}$$

$$\leq 6 \max_{I \in \mathrm{I}_C^\kappa} \xi(I) + \frac{\varepsilon}{2}.$$

Consequently,

$$\liminf_{r \to \infty} \mathbf{P}\left( \sup_{A \in \mathscr{A}} \Lambda^r(\partial_A^\kappa) \leq \varepsilon \right) \geq \liminf_{r \to \infty} \mathbf{P}\left( \left\{ \max_{I \in \mathrm{I}_C^\kappa} \Lambda^r(I) < \frac{\varepsilon}{12} \right\} \cap \Omega_0 \right).$$

Apply Fatou's lemma to the right-hand side to obtain

$$(4.42) \qquad \begin{aligned} &\liminf_{r \to \infty} \mathbf{P}\left( \sup_{A \in \mathscr{A}} \Lambda^r(\partial_A^\kappa) \leq \varepsilon \right) \\ &\qquad \geq \mathbf{P}\left( \left\{ \limsup_{r \to \infty} \max_{I \in \mathrm{I}_C^\kappa} \Lambda^r(I) < \frac{\varepsilon}{12} \right\} \cap \Omega_0 \right). \end{aligned}$$

Each $I \in \mathrm{I}_C^\kappa$ is closed, so the Portmanteau theorem implies that

$$(4.43) \qquad\qquad \limsup_{r \to \infty} \max_{I \in \mathrm{I}_C^\kappa} \Lambda^r(I) \leq \max_{I \in \mathrm{I}_C^\kappa} \Lambda(I) \qquad \text{a.s.}$$

Combining (4.42) and (4.43) yields

$$\liminf_{r \to \infty} \mathbf{P}\left( \sup_{A \in \mathscr{A}} \Lambda^r(\partial_A^\kappa) \leq \varepsilon \right) \geq \mathbf{P}\left( \left\{ \max_{I \in \mathrm{I}_C^\kappa} \Lambda(I) < \frac{\varepsilon}{12} \right\} \cap \Omega_0 \right).$$

Thus, by (4.40) and (4.41), it remains to show that

$$(4.44) \qquad\qquad \lim_{\kappa \to 0} \sup_{\xi \in \mathbf{K}_\vartheta} \max_{I \in \mathrm{I}_C^\kappa} \xi(I) = 0.$$

If (4.44) fails, then there exist $\kappa_n \to 0$, $\{\xi_n\} \subset \mathbf{K}_\vartheta$, $\{I_n : I_n \in \mathbb{I}_C^{\kappa_n}\}$ and $\delta > 0$, such that

$$\inf_n \xi_n(I_n) > \delta$$

and such that $\{I_n\}$ are either all vertical or all horizontal. Suppose that they are all vertical. For each $n$, choose $x_n \in I_n \cap \{\mathbb{R}_+ \times \{0\}\}$. Since $C$ and $\mathbf{K}_\vartheta$ are compact, assume (by passing to a subsequence if necessary) that $x_n \to x$ and $\xi_n \xrightarrow{\mathbf{w}} \xi$. For $w > 0$, let $I_x^w = [(x-w)^+, x+w] \times \mathbb{R}$. Then for each $w > 0$ and $n$ sufficiently large, $I_n \subset I_x^w$. Thus,

$$\liminf_{n\to\infty} \xi_n(I_x^w) > \delta.$$

Since $I_x^w$ is closed, the Portmanteau theorem implies that for each $w > 0$,

$$\xi(I_x^w) \ge \liminf_{n\to\infty} \xi_n(I_x^w) > \delta.$$

Letting $I_x = \bigcap_{w>0} I_x^w = \{x\} \times \mathbb{R}$, deduce that

$$(4.45) \qquad\qquad \xi(I_x) \ge \delta.$$

However, $\xi \in \mathbf{M}_\vartheta$ because $\mathbf{M}_\vartheta$ is closed (Lemma 4.5). So $\xi = \vartheta_e^z$ for some $z > 0$; by Definition 2.1,

$$(4.46) \qquad\qquad \xi(I_x) = \alpha \int_0^\infty \vartheta(\{x + uz^{-1}\} \times \mathbb{R})\, du.$$

Since $\vartheta(\cdot \times \mathbb{R})$ is a probability measure on $\mathbb{R}_+$, it has at most countably many atoms. Thus,

$$(4.47) \qquad\qquad \xi(I_x) = 0,$$

contradicting (4.45). The argument is identical if $\{I_n\}$ are all horizontal. Conclude that (4.44) must hold.  $\square$

The regularity result is now shown for the entire state descriptor $\bar{\mathcal{Z}}^r(\cdot)$. Assume (A) and let $T, L > 1$ and $\varepsilon, \eta \in (0, 1)$.

LEMMA 4.12.   *There exists $\kappa > 0$ such that*

$$(4.48) \qquad \liminf_{r\to\infty} \mathbf{P}^r\left( \sup_{A \in \mathscr{A}} \sup_{t \in [0, \lfloor rT \rfloor + L]} \bar{\mathcal{Z}}^r(t)(\partial_A^\kappa) \le \varepsilon \right) \ge 1 - \eta.$$

PROOF.   By Lemmas 4.11 and 4.10, there exist $\kappa_0 > 0$ and compact $\mathbf{K} \subset \mathbf{M}$ such that the events

$$\Omega_1^r = \left\{ \sup_{A \in \mathscr{A}} \bar{\mathcal{Z}}^r(0)(\partial_A^{\kappa_0}) \le \frac{\varepsilon}{2} \right\},$$

$$\Omega_2^r = \{ \bar{\mathcal{Z}}^r(t) \in \mathbf{K} \text{ for all } t \in [0, \lfloor rT \rfloor + L] \}$$

satisfy

$$\liminf_{r \to \infty} \mathbf{P}^r(\Omega_1^r \cap \Omega_2^r) \geq 1 - \frac{\eta}{2}. \tag{4.49}$$

Let $R_K = \{w \in \mathbb{H}_+ : \|w\| \leq K\}$ for each $K \geq 0$. Since $\mathbf{K}$ is compact, there exist $M \in (1, \infty)$ and $K < \infty$ such that

$$\sup_{\xi \in \mathbf{K}} \langle 1, \xi \rangle \leq M, \tag{4.50}$$

$$\sup_{\xi \in \mathbf{K}} \xi(R_K^c) \leq \frac{\varepsilon}{2} \tag{4.51}$$

(see, e.g., [12], Theorem A 7.5). Define $\alpha^* = \sup_{r \in \mathcal{R}} \alpha^r$, which is finite by (2.22). Let $h = \varepsilon(8\alpha^*)^{-1} \wedge (L-1)/2$, let $\kappa = \kappa_0 \wedge h(2M)^{-1}$ and let $\delta = \varepsilon((8\lceil KMh^{-1} \rceil)^{-1} \wedge 2^{-1})$. By (4.49) and Lemmas 4.7 and 4.9, there exist events

$$\Omega_3^r \subset \left\{ \sup_{A \in \mathscr{A}} \sup_{\substack{m \leq \lfloor rT \rfloor, \\ [s,t] \subset [0,L]}} |\bar{\mathcal{L}}_{s,t}^{r,m}(A) - \alpha^r(t-s)\breve{\vartheta}^r(A)| \leq \delta \right\}, \qquad r \in \mathcal{R},$$

such that the events $\Omega_0^r = \Omega_1^r \cap \Omega_2^r \cap \Omega_3^r$ satisfy

$$\liminf_{r \to \infty} \mathbf{P}^r(\Omega_0^r) \geq 1 - \eta. \tag{4.52}$$

Let $\Omega_*^r$ denote the event in (4.48). By (4.52), it suffices to show that $\Omega_0^r \subset \Omega_*^r$ for all $r \in \mathcal{R}$. Let $\omega \in \Omega_0^r$. For the remainder of the proof, all random objects are evaluated at this $\omega$.

Consider any $t \in [0, \lfloor rT \rfloor + L]$ and $A \in \mathscr{A}$. Define

$$\tau_1 = \sup\{s \leq t : \langle 1, \bar{\mathcal{Z}}^r(s) \rangle = 0\} \tag{4.53}$$

if the supremum exists and define $\tau_1 = 0$ otherwise. Let $\tau = \tau_1 \vee (t - KM)$. The first step is to show that

$$\bar{\mathcal{Z}}^r(\tau)(\partial_A^\kappa + (\bar{S}_{\tau,t}^r, t - \tau)) \leq \frac{\varepsilon}{2}. \tag{4.54}$$

If $\tau = 0$, this follows from the definition of $\Omega_1^r$ because

$$\partial_A^\kappa + (\bar{S}_{\tau,t}^r, t - \tau) \subset \partial_{A + (\bar{S}_{\tau,t}^r, t - \tau)}^\kappa$$

and because $\kappa \leq \kappa_0$ and $\mathscr{A}$ is closed under positive translation. Suppose $\tau = \tau_1 > 0$. Then $\tau = m + s$ for some $m \leq \lfloor rT \rfloor$ and $s \in (0, L]$. By (4.53), there exists a sequence $\{s_n\} \subset [0, s]$ such that $\langle 1, \bar{\mathcal{Z}}^{r,m}(s_n) \rangle = 0$ for all $n$. By (4.22) and the definition of $\Omega_3^r$,

$$\begin{aligned}
\bar{\mathcal{Z}}^r(\tau)(\partial_A^\kappa + (\bar{S}_{\tau,t}^r, t - \tau)) &= \bar{\mathcal{Z}}^{r,m}(s)(\partial_A^\kappa + (\bar{S}_{\tau,t}^r, t - \tau)) \\
&\leq \bar{\mathcal{Z}}^{r,m}(s_n)(\mathbb{H}_+) + \bar{\mathcal{L}}_{s_n,s}^{r,m}(\mathbb{H}_+) \\
&\leq \alpha^r(s - s_n) + \delta
\end{aligned}$$

for all $n \in \mathbb{N}$. Letting $s_n \uparrow s$,

$$\bar{\mathcal{Z}}^r(\tau)(\partial_A^\kappa + (\bar{S}_{\tau,t}^r, t - \tau)) \leq \delta \leq \frac{\varepsilon}{2}.$$

Suppose that $\tau = t - KM$. Since $\langle 1, \tilde{\mathcal{Z}}^r(s) \rangle > 0$ for all $s \in (\tau, t]$, the definition of $\Omega_2^r$ and (4.50) imply that

$$\bar{S}_{\tau,t}^r = \int_{t-KM}^t \langle 1, \bar{\mathcal{Z}}^r(s) \rangle^{-1} \, ds \geq K.$$

So, by definition of $\Omega_2^r$ and (4.51),

$$\bar{\mathcal{Z}}^r(\tau)(\partial_A^\kappa + (\bar{S}_{\tau,t}^r, t - \tau)) \leq \bar{\mathcal{Z}}^r(\tau)(R_K^c) \leq \frac{\varepsilon}{2}.$$

The preceding three cases prove (4.54).

Using the dynamic equation (4.20), we obtain

$$(4.55) \qquad \begin{aligned} \tilde{\mathcal{Z}}^r(t)(\partial_A^\kappa) &= \tilde{\mathcal{Z}}^r(\tau)(\partial_A^\kappa + (\bar{S}_{\tau,t}^r, t - \tau)) \\ &\quad + \frac{1}{r} \sum_{i=r\bar{E}^r(\tau)+1}^{r\bar{E}^r(t)} 1_{\partial_A^\kappa}^+(\bar{v}_i^r(t), \bar{l}_i^r(t)). \end{aligned}$$

Let $I$ denote the second right-hand term in (4.55). By (4.54), it remains to show that $I \leq \varepsilon/2$. Let $N = \lceil (t - \tau)h^{-1} \rceil$ and for each $n \in \{0, \ldots, N-1\}$, let $t_n = \tau + nh$ and $t^n = t_{n+1} \wedge t$. Then, using the inequality $1_{\partial_A^\kappa}^+(\cdot, \cdot) \leq 1_{\partial_A^\kappa}(\cdot, \cdot)$,

$$(4.56) \qquad I \leq \sum_{n=0}^{N-1} \frac{1}{r} \sum_{i=r\bar{E}^r(t_n)+1}^{r\bar{E}^r(t^n)} 1_{\partial_A^\kappa}(\bar{v}_i^r(t), \bar{l}_i^r(t)).$$

Consider $n \in \{0, \ldots, N-1\}$ and $i$ such that $U_i^r r^{-1} \in (t_n, t^n]$. Observe that

$$(4.57) \qquad \bar{S}_{t^n,t}^r \leq \bar{S}_{U_i^r r^{-1}, t}^r \leq \bar{S}_{t_n,t}^r.$$

By (2.3), (2.5), (2.14) and (2.15), we have

$$(4.58) \qquad 1_{\partial_A^\kappa}(\bar{v}_i^r(t), \bar{l}_i^r(t)) = 1_{\partial_A^\kappa + (\bar{S}_{U_i^r r^{-1}, t}^r, t - U_i^r r^{-1})}(v_i^r, l_i^r r^{-1}).$$

So, letting

$$A_n^- = A + (\bar{S}_{t^n,t}^r - \kappa, t - t^n - \kappa) \cap \mathbb{H}_+,$$

$$A_n^+ = A + (\bar{S}_{t_n,t}^r + \kappa, t - t_n + \kappa),$$

$$A_n = A_n^- \setminus A_n^+,$$

it follows from (4.57) and (4.58) that

$$(4.59) \qquad 1_{\partial_A^\kappa}(\bar{v}_i^r(t), \bar{l}_i^r(t)) \leq 1_{A_n}(v_i^r, l_i^r r^{-1}).$$

Conclude from (4.56) and (4.59) that

$$I \leq \sum_{n=0}^{N-1} \frac{1}{r} \sum_{i=r\bar{E}^r(t_n)+1}^{r\bar{E}^r(t^n)} 1_{A_n}(v_i^r, l_i^r r^{-1}) = \sum_{n=0}^{N-1} (\bar{\mathcal{L}}_{t_n, t^n}^r(A_n^-) - \bar{\mathcal{L}}_{t_n, t^n}^r(A_n^+)).$$

Note that $t^n - t_n \leq h < L$ for all $n$, so $\bar{\mathcal{L}}_{t_n, t^n}^r$ can be rewritten $\bar{\mathcal{L}}_{s,t}^{r,m}$ for some $m \leq \lfloor rT \rfloor$ and $[s,t] \subset [0, L]$. Since $A_n^-, A_n^+ \in \mathscr{A}$ for all $n < N$, we can deduce from the definition of $\Omega_3^r$ that

$$I \leq \sum_{n=0}^{N-1} (\alpha^r h \breve{\vartheta}^r(A_n) + 2\delta).$$

By the definitions of $\alpha^*$ and $N$ and the fact that $t - \tau \leq KM$,

$$I \leq \alpha^* h \sum_{n=0}^{N-1} \breve{\vartheta}^r(A_n) + \lceil KMh^{-1} \rceil 2\delta.$$

This implies, by choice of $\delta$, that

(4.60) $$I \leq \alpha^* h \sum_{n=0}^{N-1} \breve{\vartheta}^r(A_n) + \frac{\varepsilon}{4}.$$

If $n \in \{0, \dots, N-3\}$, then

$$\bar{S}_{t_{n+1}, t_{n+2}}^r \geq hM^{-1} \geq 2\kappa$$

because $0 < \langle 1, \bar{\mathcal{Z}}^r(s) \rangle \leq M$ for all $s \in (\tau, t]$ and because $h \geq \kappa 2M$, by definition of $\kappa$. Thus, for all $n \in \{0, \dots, N-3\}$,

$$\bar{S}_{t^n, t}^r - \kappa = \bar{S}_{t_{n+1}, t_{n+2}}^r + \bar{S}_{t_{n+2}, t}^r - \kappa \geq \bar{S}_{t_{n+2}, t}^r + \kappa.$$

Deduce that $A_n^- \subset A_{n+2}^+$ for all $n \in \{0, \dots, N-3\}$ and, consequently, that $A_n \cap A_{n+2} = \varnothing$. Thus, since $\breve{\vartheta}^r$ is a probability measure, $\sum_{n=0}^{\lfloor (N-1)/2 \rfloor} \breve{\vartheta}^r(A_{2n})$ and $\sum_{n=0}^{\lfloor (N-2)/2 \rfloor} \breve{\vartheta}^r(A_{2n+1})$ are both bounded by one. Conclude from (4.60) that

$$I \leq 2\alpha^* h + \frac{\varepsilon}{4},$$

which implies, by choice of $h$, that $I \leq \varepsilon/2$.   $\square$

4.6. *Oscillation bound.* This section contains the main oscillation bound that constitutes the first of two main ingredients needed to prove tightness of the state descriptors.

DEFINITION 4.13.   For each $L > 1$, $\zeta(\cdot) \in \mathbf{D}([0,\infty), \mathbf{M})$ and $\delta > 0$, define the modulus of continuity of $\zeta(\cdot)$ on $[0, L]$ as

$$(4.61) \qquad \mathbf{w}_L(\zeta(\cdot), \delta) = \sup_{t \in [0, L-\delta]} \sup_{h \in [0, \delta]} \mathbf{d}[\zeta(t+h), \zeta(t)].$$

Denote the modulus of continuity of $\zeta(\cdot)$ on $[0, L)$ by

$$(4.62) \qquad \mathbf{w}_{L-}(\zeta(\cdot), \delta) = \sup_{t \in [0, L-\delta)} \sup_{h \in [0, \delta]} \mathbf{d}[\zeta(t+h), \zeta(t)].$$

Assume (A) and let $T, L > 1$ and $\varepsilon, \eta \in (0, 1)$.

LEMMA 4.14.   *There exists $\delta > 0$ such that*

$$(4.63) \qquad \liminf_{r \to \infty} \mathbf{P}^r \left( \max_{m \leq \lfloor rT \rfloor} \mathbf{w}_L(\bar{\mathcal{Z}}^{r,m}(\cdot), \delta) \leq \varepsilon \right) \geq 1 - \eta.$$

PROOF.   By Lemmas 4.12 and 4.7, there exists $\kappa \in (0, 1)$ such that for each fixed $\delta > 0$, the events

$$\Omega_1^r = \left\{ \max_{m \leq \lfloor rT \rfloor} \sup_{t \in [0, L]} \bar{\mathcal{Z}}^{r,m}(t)([0, \kappa] \times \mathbb{R}) \leq \frac{\varepsilon}{4} \right\},$$

$$\Omega_2^r = \left\{ \max_{m \leq \lfloor rT \rfloor} \sup_{t \in [0, L-\delta]} \bar{\mathcal{L}}_{t, t+\delta}^{r,m}(\mathbb{H}_+) \leq 2\alpha\delta \right\},$$

$$\Omega_0^r = \Omega_1^r \cap \Omega_2^r$$

satisfy

$$(4.64) \qquad \liminf_{r \to \infty} \mathbf{P}^r(\Omega_0^r) \geq 1 - \eta.$$

Fix $\delta = \kappa \varepsilon^2 (8(\alpha \vee 1))^{-1}$ and let $\Omega_*^r$ be the event in (4.63). By (4.64), it suffices to show that $\Omega_0^r \subset \Omega_*^r$ for each $r$. Fix $r \in \mathcal{R}$ and $\omega \in \Omega_0^r$. For the remainder of the proof, all random objects are evaluated at this $\omega$. Fix $m \leq \lfloor rT \rfloor$, $t \in [0, L-\delta]$, $h \in [0, \delta]$ and let $B \subset \mathbb{H}_+$ be closed. By definition of the Prohorov metric $\mathbf{d}[\cdot, \cdot]$, it suffices to show the two inequalities

$$(4.65) \qquad \bar{\mathcal{Z}}^{r,m}(t)(B) \leq \bar{\mathcal{Z}}^{r,m}(t+h)(B^\varepsilon) + \varepsilon,$$

$$(4.66) \qquad \bar{\mathcal{Z}}^{r,m}(t+h)(B) \leq \bar{\mathcal{Z}}^{r,m}(t)(B^\varepsilon) + \varepsilon.$$

To show (4.65), we use the definition of $\Omega_1^r$ to write

$$(4.67) \qquad \begin{aligned} \bar{\mathcal{Z}}^{r,m}(t)(B) &\leq \bar{\mathcal{Z}}^{r,m}(t)([0, \kappa] \times \mathbb{R}) + \bar{\mathcal{Z}}^{r,m}(t)(B \cap ((\kappa, \infty) \times \mathbb{R})) \\ &\leq \frac{\varepsilon}{4} + \bar{\mathcal{Z}}^{r,m}(t)(B \cap ((\kappa, \infty) \times \mathbb{R})). \end{aligned}$$

Let $I = \{s \in [t, t+h] : \langle 1, \bar{\mathcal{Z}}^{r,m}(s) \rangle < \varepsilon/4\}$. Suppose that $I = \varnothing$. Then $\langle 1, \bar{\mathcal{Z}}^{r,m}(s) \rangle \geq \varepsilon/4$ for all $s \in [t, t+h]$, which implies that

$$(4.68) \quad \|(\bar{S}_{t,t+h}^{r,m}, h)\| \leq \int_t^{t+\delta} \langle 1, \bar{\mathcal{Z}}^{r,m}(s) \rangle^{-1}\, ds + \delta \leq \frac{4\delta}{\varepsilon} + \delta < \varepsilon \wedge \kappa.$$

Consequently, $w \in B \cap ((\kappa, \infty) \times \mathbb{R})$ implies that $w - (\bar{S}_{t,t+h}^{r,m}, h) \in B^\varepsilon$ and so

$$(4.69) \quad B \cap ((\kappa, \infty) \times \mathbb{R}) \subset B^\varepsilon + (\bar{S}_{t,t+h}^{r,m}, h).$$

Deduce from (4.67) that

$$\bar{\mathcal{Z}}^{r,m}(t)(B) \leq \frac{\varepsilon}{4} + \bar{\mathcal{Z}}^{r,m}(t)(B^\varepsilon + (\bar{S}_{t,t+h}^{r,m}, h)).$$

Applying (4.23), we obtain

$$(4.70) \quad \bar{\mathcal{Z}}^{r,m}(t)(B) \leq \frac{\varepsilon}{4} + \bar{\mathcal{Z}}^{r,m}(t+h)(B^\varepsilon).$$

Suppose that $I \neq \varnothing$ and let $\tau = \inf I$. Then $\langle 1, \bar{\mathcal{Z}}^{r,m}(\tau) \rangle \leq \varepsilon/4$, by right continuity. If $\tau > t$, then $\langle 1, \bar{\mathcal{Z}}^{r,m}(s) \rangle \geq \varepsilon/4$ for all $s \in [t, \tau)$, so

$$(4.71) \quad \bar{S}_{t,\tau}^{r,m} = \int_t^\tau \langle 1, \bar{\mathcal{Z}}^{r,m}(s) \rangle^{-1}\, ds \leq \frac{4(\tau - t)}{\varepsilon} \leq \frac{4\delta}{\varepsilon} < \kappa.$$

By (4.67) and (4.71),

$$\bar{\mathcal{Z}}^{r,m}(t)(B) \leq \frac{\varepsilon}{4} + \bar{\mathcal{Z}}^{r,m}(t)((\kappa, \infty) \times \mathbb{R})$$

$$\leq \frac{\varepsilon}{4} + \bar{\mathcal{Z}}^{r,m}(t)((\bar{S}_{t,\tau}^{r,m}, \infty) \times \mathbb{R}).$$

Applying (4.23), we obtain

$$(4.72) \quad \bar{\mathcal{Z}}^{r,m}(t)(B) \leq \frac{\varepsilon}{4} + \bar{\mathcal{Z}}^{r,m}(\tau)(\mathbb{H}_+) \leq \frac{\varepsilon}{2}.$$

Therefore, (4.65) follows because either (4.70) or (4.72) holds.

To show (4.66), we use (4.22) and the definitions of $\Omega_2^r$ and $\delta$ to obtain

$$(4.73) \quad \begin{aligned} \bar{\mathcal{Z}}^{r,m}(t+h)(B) &\leq \bar{\mathcal{Z}}^{r,m}(t)(B + (\bar{S}_{t,t+h}^{r,m}, h)) + \bar{\mathcal{L}}_{t,t+h}^{r,m}(\mathbb{H}_+) \\ &\leq \bar{\mathcal{Z}}^{r,m}(t)(B + (\bar{S}_{t,t+h}^{r,m}, h)) + \frac{\varepsilon}{4}. \end{aligned}$$

If $I = \varnothing$, then (4.68) implies that $B + (\bar{S}_{t,t+h}^{r,m}, h) \subset B^\varepsilon$. Therefore, (4.73) yields

$$\bar{\mathcal{Z}}^{r,m}(t+h)(B) \leq \bar{\mathcal{Z}}^{r,m}(t)(B^\varepsilon) + \frac{\varepsilon}{4}.$$

If $I \neq \varnothing$, then by (4.22) and the definitions of $\Omega_2^r$ and $\delta$,

$$\bar{\mathcal{Z}}^{r,m}(t+h)(B) \leq \bar{\mathcal{Z}}^{r,m}(\tau)(\mathbb{H}_+) + \bar{\mathcal{L}}_{\tau,t+h}^{r,m}(\mathbb{H}_+) \leq \frac{\varepsilon}{4} + 2\alpha\delta \leq \frac{\varepsilon}{2}.$$

In both cases, (4.66) holds.

Conclude from (4.65) and (4.66) that

$$\mathbf{d}[\bar{\mathcal{Z}}^{r,m}(t), \bar{\mathcal{Z}}^{r,m}(t+h)] \leq \varepsilon.$$

Since $m \leq \lfloor rT \rfloor$, $t \in [0, L-\delta]$ and $h \in [0, \delta]$ were arbitrary,

$$\max_{m \leq \lfloor rT \rfloor} \mathbf{w}_L(\bar{\mathcal{Z}}^{r,m}(\cdot), \delta) \leq \varepsilon,$$

which implies that $\omega \in \Omega_*^r$.   □

4.7. *Precompactness.* This section introduces a sequence of events on which the shifted fluid scaled state descriptors $\{\bar{\mathcal{Z}}^{r,m}(\cdot) : m \leq \lfloor rT \rfloor\}$ have a desired list of properties. It will be shown that the tail of this sequence has arbitrarily high probability and that certain sequences of sample paths chosen from the events are precompact. This will be used in the next section to construct fluid approximations to the shifted fluid scaled state descriptors.

Assume (A) and let $T > 1$. Fix a positive constant $M$ and a compact set $\mathbf{K} \subset \mathbf{M}$. For each $L > 1$, let $\mathcal{C}_L = (\{\kappa_j\}_{j=1}^{\infty}, \{\delta_k\}_{k=1}^{\infty})$ be a collection, depending on $L$, where $\kappa_j, \delta_k$ are positive constants with $\delta_k \to 0$ as $k \to \infty$.

DEFINITION 4.15. For each $L > 1$ and each $r \in \mathcal{R}$ and $n \in \mathbb{N}$, define events

$$\Omega^{r,1} = \left\{ \max_{m \leq \lfloor rT \rfloor} \sup_{t \in [0,L]} \langle 1 \vee \chi, \bar{\mathcal{Z}}^{r,m}(t) \rangle \leq M \right\},$$

$$\Omega^{r,2} = \{ \bar{\mathcal{Z}}^{r,m}(t) \in \mathbf{K} \quad \text{for all } m \leq \lfloor rT \rfloor \text{ and } t \in [0,L] \},$$

$$\Omega_n^{r,3} = \left\{ \max_{m \leq \lfloor rT \rfloor} \sup_{t \in [0,L]} |\langle 1, \bar{\mathcal{Z}}^{r,m}(t) \rangle - \langle 1, \bar{\mathcal{Z}}^{r,m}(0) \rangle| \leq \frac{1}{n} \right\},$$

$$\Omega_n^{r,4} = \bigcap_{j=1}^{n} \left\{ \sup_{A \in \mathscr{A}} \max_{m \leq \lfloor rT \rfloor} \sup_{t \in [0,L]} \bar{\mathcal{Z}}^{r,m}(t)(\partial_A^{\kappa_j}) \leq \frac{1}{j} \right\},$$

$$\Omega_n^{r,5} = \bigcap_{k=1}^{n} \left\{ \max_{m \leq \lfloor rT \rfloor} \mathbf{w}_L(\bar{\mathcal{Z}}^{r,m}(\cdot), \delta_k) \leq \frac{1}{k} \right\},$$

$$\Omega_n^{r,6} \subset \left\{ \sup_{A \in \mathscr{A}} \sup_{\substack{m \leq \lfloor rT \rfloor, \\ [s,t] \subset [0,L]}} |\bar{\mathcal{L}}_{s,t}^{r,m}(A) - \alpha^r(t-s)\vartheta^r(A)| \leq \frac{1}{n} \right\},$$

$$\Omega_{n,L}^r(M, \mathbf{K}, \mathcal{C}_L) = \Omega^{r,1} \cap \Omega^{r,2} \cap \Omega_n^{r,3} \cap \Omega_n^{r,4} \cap \Omega_n^{r,5} \cap \Omega_n^{r,6},$$

where for each fixed $n \in \mathbb{N}$, the events $\{\Omega_n^{r,6} : r \in \mathcal{R}\}$ are chosen so that $\lim_{r \to \infty} \mathbf{P}^r(\Omega_n^{r,6}) = 1$ (see Lemmas 4.7 and 4.9).

Fix $M$, $\mathbf{K}$ and collections $\{\mathcal{C}_L : L > 1\}$. Then for each $L > 1$, the events $\{\Omega_{n,L}^r(M, \mathbf{K}, \mathcal{C}_L) : r \in \mathcal{R}, \, n \in \mathbb{N}\}$ form an array indexed by $\mathcal{R} \times \mathbb{N}$. The next lemma asserts that $M$, $\mathbf{K}$ and $\{\mathcal{C}_L : L > 1\}$ can be chosen so that for each $L > 1$, the tail of every $\mathcal{R}$-indexed column in this array has arbitrarily high probability.

Assume (A) and let $T > 1$ and $0 < \eta < 1$.

LEMMA 4.16. *There exist $M^* > 0$, a compact $\mathbf{K}^* \subset \mathbf{M}$ and collections $\{\mathcal{C}_L^* : L > 1\}$ such that for each $L > 1$ and each fixed $n \in \mathbb{N}$,*

$$(4.74) \qquad \liminf_{r \to \infty} \mathbf{P}^r(\Omega_{n,L}^r(M^*, \mathbf{K}^*, \mathcal{C}_L^*)) \geq 1 - \frac{\eta}{2}.$$

PROOF. Use Corollary 4.3 and Lemma 4.10 to choose $M^*$ and $\mathbf{K}^*$ so that for all $L > 1$ and each fixed $n$,

$$(4.75) \qquad \liminf_{r \to \infty} \mathbf{P}^r(\Omega^{r,1} \cap \Omega^{r,2} \cap \Omega_n^{r,3}) \geq 1 - \frac{\eta}{4}.$$

Let $L > 1$. By Lemmas 4.12, 4.14, 4.7 and 4.9, there exists $\mathcal{C}_L^* = (\{\kappa_j\}, \{\delta_k\})$ such that for each fixed $n$,

$$(4.76) \qquad \liminf_{r \to \infty} \mathbf{P}^r(\Omega_n^{r,4} \cap \Omega_n^{r,5} \cap \Omega_n^{r,6}) \geq 1 - \frac{\eta}{4}. \qquad \square$$

The preceding lemma guarantees that for each $L > 1$, there exists a "diagonal" sequence of events, indexed by $r \in \mathcal{R}$, such that the tail of this sequence has arbitrarily high probability.

Assume (A) and let $T > 1$ and $\eta \in (0, 1)$. Using Lemma 4.16 choose a constant $M^*$, a compact set $\mathbf{K}^*$ and collections $\{\mathcal{C}_L^* : L > 1\}$ so that (4.74) holds for each $L$ and $n$.

DEFINITION 4.17. For each $L > 1$ and $r \in \mathcal{R}$, let

$$\mathcal{N}_L^r = \{n \in \mathbb{N} : \mathbf{P}^r(\Omega_{n,L}^r(M^*, \mathbf{K}^*, \mathcal{C}_L^*)) \geq 1 - \eta\}$$

and let

$$(4.77) \qquad n(L, r) = \begin{cases} \sup \mathcal{N}_L^r \wedge \lfloor r \rfloor, & \mathcal{N}_L^r \neq \varnothing, \\ 1, & \mathcal{N}_L^r = \varnothing. \end{cases}$$

For each $L > 1$ and $r \in \mathcal{R}$, define

$$(4.78) \qquad \Omega_L^r = \Omega_{n(L,r),L}^r(M^*, \mathbf{K}^*, \mathcal{C}_L^*).$$

Note that in the above definition, $\sup \mathcal{N}_L^r$ could equal infinity for some $r$. In this case, it is sufficient to use $\Omega_L^r = \Omega_{\lfloor r \rfloor, L}^r(M^*, \mathbf{K}^*, \mathcal{C}_L^*)$, which is the reason for including the minimum with $\lfloor r \rfloor$ in the definition. This does not affect the property that $n(L, r) \in \mathcal{N}_L^r$ for $\mathcal{N}_L^r \neq \varnothing$: if $m, n \in \mathbb{N}$ with $m \leq n$, then $\Omega_{n,L}^r(M^*, \mathbf{K}^*, C_L^*) \subset \Omega_{m,L}^r(M^*, \mathbf{K}^*, C_L^*)$. So, for $m \leq n$, $n \in \mathcal{N}_L^r$ implies that $m \in \mathcal{N}_L^r$. In particular, $\mathcal{N}_L^r \neq \varnothing$ implies that $n(L, r) \in \mathcal{N}_L^r$.

Also, note that $\Omega_L^r$ may be empty for some $r$, possibly when $\mathcal{N}_L^r = \varnothing$. This is of no concern; the following lemma implies that this is the case for at most finitely many $r$.

Assume (A) and let $T > 1$ and $\eta \in (0, 1)$.

LEMMA 4.18.   *For each $L > 1$,*

$$(4.79) \qquad\qquad \liminf_{r \to \infty} \mathbf{P}^r(\Omega_L^r) \geq 1 - \eta.$$

PROOF.   Fix $L > 1$. If $n(r, L) > 1$ for some $r \in \mathcal{R}$, then $\mathbf{P}^r(\Omega_L^r) \geq 1 - \eta$ by definition. Thus, it suffices to show that

$$(4.80) \qquad\qquad n(r, L) \to \infty \qquad \text{as } r \to \infty.$$

This follows by Definition 4.17 and Lemma 4.16.   $\square$

Using the events $\Omega_L^r$, define, for each $L > 1$, a sequence of subsets of $\mathbf{D}([0, \infty), \mathbf{M})$ that are sample paths on $[0, L)$ of the shifted, fluid scaled state descriptors $\{\bar{\mathcal{Z}}^{r,m}(\cdot) : m = 1, \ldots, \lfloor rT \rfloor\}$. Assume (A) and let $T > 1$ and $\eta \in (0, 1)$.

DEFINITION 4.19.   For each $L > 1$ and $r \in \mathcal{R}$, let $\mathscr{D}_L^r$ be the set of all $\zeta(\cdot) \in \mathbf{D}([0, \infty), \mathbf{M})$ such that for some $\omega \in \Omega_L^r$ and some $m \leq \lfloor rT \rfloor$,

$$(4.81) \qquad\qquad \zeta(t) = \begin{cases} \bar{\mathcal{Z}}_\omega^{r,m}(t), & t \in [0, L), \\ \mathbf{0}, & t \in [L, \infty). \end{cases}$$

By Lemma 4.18, $\mathscr{D}_L^r$ is empty for at most finitely many $r$.

Assume (A) and let $T > 1$ and $\eta \in (0, 1)$. For each $L > 1$, let $\{\mathscr{D}_L^r : r \in \mathcal{R}\}$ be the sequence of sets defined above. Fix $L > 1$ and suppose that $\tilde{\mathcal{R}} \subset \mathcal{R}$ is a subsequence and $\{\zeta^{\tilde{r}}(\cdot) : \tilde{r} \in \tilde{\mathcal{R}}\} \subset \mathbf{D}([0, \infty), \mathbf{M})$ is a sequence of paths such that $\zeta^{\tilde{r}}(\cdot) \in \mathscr{D}_L^{\tilde{r}}$ for each $\tilde{r} \in \tilde{\mathcal{R}}$.

THEOREM 4.20.   *The sequence $\{\zeta^{\tilde{r}}(\cdot) : \tilde{r} \in \tilde{\mathcal{R}}\}$ is precompact in $\mathbf{D}([0, \infty),$ $\mathbf{M})$.*

PROOF.   Let $\mathbf{K}^*$ and $\mathcal{C}_L^* = (\{\kappa_j\}, \{\delta_k\})$ be given by Lemma 4.16. By definition of $\mathscr{D}_L^{\tilde{r}}$ and $\Omega_L^{\tilde{r}}$ ($\Omega^{\tilde{r},2}$ in particular), $\zeta^{\tilde{r}}(t) \in \mathbf{K}^* \cup \{\mathbf{0}\}$ for all $t \geq 0$ and $\tilde{r} \in \tilde{\mathcal{R}}$. Let $k \in \mathbb{N}$. By (4.80), there exists $r_k$ such that $\tilde{r} > r_k$ implies that $n(\tilde{r}, L) > k$. For $R > 1$, let $\mathbf{w}_R'(\cdot, \cdot)$ be the modulus of continuity on $[0, R]$ used in Theorem 3.6.3 of [7]. Since $\lim_{\delta \to 0} \mathbf{w}'_R(\xi(\cdot), \delta) = 0$ for each fixed $\xi(\cdot) \in \mathbf{D}([0, \infty), \mathbf{M})$, there exists $\delta'_k \in (0, \delta_k)$ such that

$$(4.82) \qquad \max_{\tilde{r} \leq r_k} \mathbf{w}'_R(\zeta^{\tilde{r}}(\cdot), \delta'_k) < \frac{1}{k}.$$

If $\tilde{r} > r_k$, then $\Omega^{\tilde{r}, 5}_{n(\tilde{r}, L)} \subset \Omega^{\tilde{r}, 5}_k$ (see Definition 4.15). By definition of $\mathscr{D}^{\tilde{r}}_L$ and $\Omega^{\tilde{r}}_L$ ($\Omega^{\tilde{r}, 5}_{n(\tilde{r}, L)}$ in particular), we have

$$(4.83) \qquad \sup_{\tilde{r} > r_k} \mathbf{w}_{L-}(\zeta^{\tilde{r}}(\cdot), \delta'_k) \leq \sup_{\tilde{r} > r_k} \mathbf{w}_{L-}(\zeta^{\tilde{r}}(\cdot), \delta_k) \leq \frac{1}{k}.$$

Since $\mathbf{w}'_L(\zeta^{\tilde{r}}(\cdot), \delta'_k)$ is bounded above by $\mathbf{w}_{L-}(\zeta^{\tilde{r}}(\cdot), \delta'_k)$ for each $\tilde{r} \in \tilde{\mathcal{R}}$ and $k \in \mathbb{N}$ and since $\zeta^{\tilde{r}}(\cdot) \equiv \mathbf{0}$ on $[L, \infty)$ for each $\tilde{r}$, (4.83) implies that for $R > 1$,

$$(4.84) \qquad \lim_{k \to \infty} \sup_{\tilde{r} > r_k} \mathbf{w}'_R(\zeta^{\tilde{r}}(\cdot), \delta'_k) \leq \lim_{k \to \infty} \sup_{\tilde{r} > r_k} \mathbf{w}_{L-}(\zeta^{\tilde{r}}(\cdot), \delta'_k) = 0.$$

Deduce from (4.82) and (4.84) that

$$(4.85) \qquad \lim_{k \to \infty} \sup_{\tilde{r} \in \tilde{\mathcal{R}}} \mathbf{w}'_R(\zeta^{\tilde{r}}(\cdot), \delta'_k) = 0.$$

So $\{\zeta^{\tilde{r}}(\cdot) : \tilde{r} \in \tilde{\mathcal{R}}\}$ is precompact by Theorem 3.6.3 of [7]. $\quad \square$

## 5. Local fluid limits.

This section provides the remaining ingredients needed to prove Theorem 2.2. It is divided into three parts. Section 5.1 below identifies the limit set $\mathscr{D}_L$ of the sequence of sets $\{\mathscr{D}^r_L : r \in \mathcal{R}\}$ introduced in the previous section and derives several properties of the elements of this set. Section 5.2 establishes that the set $\mathscr{D}_L$ approximates the sample paths of $\tilde{\mathcal{Z}}^{r,m}(\cdot)$ on $[0, L)$ with high probability. Section 5.3 contains a result on the steady state behavior of the elements of $\mathscr{D}_L$, an important ingredient for proving state space collapse of the diffusion scaled state descriptors $\{\hat{\mathcal{Z}}^r(\cdot)\}$. The results of this section are combined in Section 6 to show state space collapse, which quickly leads to the proof of Theorem 2.2.

### 5.1. Properties.

Assume (A) and let $T > 1$ and $\eta \in (0, 1)$.

DEFINITION 5.1. For each $L > 1$, let $\mathscr{D}_L$ be the set of all $\zeta(\cdot) \in \mathbf{D}([0, \infty), \mathbf{M})$ such that there exists a subsequence $\tilde{\mathcal{R}} \subset \mathcal{R}$ and a sequence $\{\zeta^{\tilde{r}}(\cdot) : \tilde{r} \in \tilde{\mathcal{R}}\} \subset \mathbf{D}([0, \infty), \mathbf{M})$, satisfying $\zeta^{\tilde{r}}(\cdot) \in \mathscr{D}^{\tilde{r}}_L$ for each $\tilde{r} \in \tilde{\mathcal{R}}$ and

$$(5.1) \qquad \zeta^{\tilde{r}}(\cdot) \xrightarrow{J_1} \zeta(\cdot) \qquad \text{as } \tilde{r} \to \infty.$$

The set $\mathscr{D}_L$ is nonempty by Lemma 4.18 and Theorem 4.20. Call the elements of $\mathscr{D}_L$ *local fluid limits on* $[0, L)$ of the state descriptor $\tilde{\mathcal{Z}}^r(\cdot)$; the term indicates that each element $\zeta(\cdot) \in \mathscr{D}_L$ approximates some section of a sample path of $\bar{\mathcal{Z}}^r(\cdot): [0, rT] \to \mathbf{M}$ over a finite time interval of length $L$. The precise statement of this approximation appears in Section 5.2 below; the purpose of this section is to establish several properties of the elements of $\mathscr{D}_L$.

Let $\mathscr{H}$ denote the Borel subsets of $\mathbb{H}_+$. Assume (A) and let $T > 1$ and $\eta \in (0, 1)$.

THEOREM 5.2. *There exist a positive constant $M$ and a compact set $\mathbf{K} \subset \mathbf{M}$ such that for each $L > 1$ and each $\zeta(\cdot) \in \mathscr{D}_L$:*

(i) *$\zeta(t) \in \mathbf{K}$ for all $t \geq 0$;*
(ii) *$\zeta(\cdot) \equiv \mathbf{0}$ on $[L, \infty)$;*
(iii) *$\zeta(\cdot)$ is continuous on $[0, L)$;*
(iv) *$\langle 1, \zeta(t) \rangle \leq M$ for all $t \in [0, L)$;*
(v) *$\langle 1, \zeta(\cdot) \rangle$ is constant on $[0, L)$;*
(vi) *if $z = \langle 1, \zeta(0) \rangle$ is positive, then for all $t \in [0, L)$ and $B \in \mathscr{H}$,*

$$(5.2) \qquad \zeta(t)(B) = \zeta(0)(B + (tz^{-1}, t)) + \alpha \int_0^t \vartheta(B + (sz^{-1}, s)) \, ds.$$

PROOF. Let $M^*$, $\mathbf{K}^*$ and $\mathcal{C}_L^* = (\{\kappa_j\}, \{\delta_n\})$ be the objects used in Definitions 4.17 and 4.19 to define $\mathscr{D}_L^r$. Let $M = M^*$ and $\mathbf{K} = \mathbf{K}^*$. Fix $L > 1$ and $\zeta(\cdot) \in \mathscr{D}_L$. By Definitions 4.19 and 5.1, there exist a subsequence $\tilde{\mathcal{R}} \subset \mathcal{R}$ and a sequence $\{\zeta^{\tilde{r}}(\cdot) : \tilde{r} \in \tilde{\mathcal{R}}\}$ such that $\zeta^{\tilde{r}}(\cdot) \in \mathscr{D}_L^{\tilde{r}}$ for each $\tilde{r}$ and such that (5.1) holds. To ease notation, assume without loss of generality that $\tilde{\mathcal{R}} = \mathcal{R}$. Properties (i) and (ii) follow immediately from Definitions 4.17, 4.19 and 5.1. Property (iii) follows from the proof of Theorem 4.20, in particular, from (4.83).

Since $\zeta(\cdot)$ is continuous on $[0, L)$, (5.1) implies that for all $t \in [0, L)$,

$$(5.3) \qquad \qquad \zeta^r(t) \xrightarrow{\mathbf{w}} \zeta(t) \qquad \text{as } r \to \infty.$$

Consequently, for all $t \in [0, L)$,

$$(5.4) \qquad \qquad \langle 1, \zeta(t) \rangle = \lim_{r \to \infty} \langle 1, \zeta^r(t) \rangle \leq M,$$

by the definition of $\mathscr{D}_L^r$ ($\Omega^{r,1}$ in particular). This proves (iv). Property (v) follows from (5.3) and the definition of $\mathscr{D}_L^r$ ($\Omega_n^{r,3}$ in particular).

To prove (vi), suppose that $z = \langle 1, \zeta(0) \rangle$ is positive and fix $t \in [0, L)$ and $A \in \mathscr{A}$. An extension to Borel sets $B \in \mathscr{H}$ will be established at the end of the proof. For each $r \in \mathcal{R}$, there exists $\omega \in \Omega_L^r$ and $m \leq \lfloor rT \rfloor$ such that, evaluated at $\omega$,

$$\zeta^r(t)(A) = \zeta^r(0)(A + (\bar{S}^{r,m}(t), t))$$

$$(5.5)$$

$$+ \frac{1}{r} \sum_{i=r\bar{E}^{r,m}(0)+1}^{r\bar{E}^{r,m}(t)} 1_A^+ (\bar{v}_i^{r,m}(t), \bar{l}_i^{r,m}(t));$$

see Definition 4.19 and the dynamic equation (4.21). By (4.83) and (5.3),

$$(5.6) \qquad \lim_{r\to\infty} \| \langle 1, \zeta^r(\cdot) \rangle - \langle 1, \zeta(\cdot) \rangle \|_L = 0.$$

By property (v), $\langle 1, \zeta(u) \rangle = z > 0$ for all $u \in [0, t]$. Therefore, (5.6) and the bounded convergence theorem imply that for all $s \in [0, t]$,

$$(5.7) \qquad \begin{aligned} \lim_{r\to\infty} \bar{S}^{r,m}(s) &= \lim_{r\to\infty} \int_0^s \varphi(\langle 1, \zeta^r(u) \rangle) \, du \\ &= \int_0^s \langle 1, \zeta(u) \rangle^{-1} \, du \\ &= sz^{-1}. \end{aligned}$$

For $s \in [0, t]$, define $S_{s,t} = (t - s)z^{-1}$. Since $s \mapsto sz^{-1}$ is continuous on $[0, t]$, the convergence in (5.7) is uniform on $[0, t]$. For $\delta > 0$, there exists $r_\delta \in \mathcal{R}$ such that

$$(5.8) \qquad \sup_{s\in[0,t]} |\bar{S}_{s,t}^{r,m} - S_{s,t}| \leq \delta \qquad \text{for all } r > r_\delta.$$

The following lemma will be needed.

LEMMA 5.3.   *For all $A \in \mathscr{A}$ and $s \in [0, L)$,*

$$(5.9) \qquad \zeta(s)(\partial_A) = 0.$$

PROOF.   Fix $A \in \mathscr{A}$ and observe that $\partial_A \subset \partial_A^{\kappa_n}$ for all $n \in \mathbb{N}$. For all $t \in [0, L)$ and $n \in \mathbb{N}$, (5.3), the Portmanteau theorem and the definition of $\mathscr{D}_L^r$ ($\Omega_n^{r,4}$ in particular) yield

$$\zeta(t)(\partial_A) \leq \zeta(t)(\partial_A^{\kappa_n}) \leq \liminf_{r\to\infty} \zeta^r(t)(\partial_A^{\kappa_n}) \leq \frac{1}{n}.$$

Letting $n \to \infty$ proves (5.9).   $\square$

Continuing with the proof of Theorem 5.2(vi), we use (5.7) and the definition of $\mathscr{D}_L^r$ ($\Omega_n^{r,4}$ in particular) to obtain

$$(5.10) \qquad |\zeta^r(0)(A + (\bar{S}^{r,m}(t), t)) - \zeta^r(0)(A + (tz^{-1}, t))| \to 0$$

as $r \to \infty$. Deduce from (5.3), (5.10), Lemma 5.3 and the Portmanteau theorem that as $r \to \infty$,

$$(5.11) \quad \zeta^r(t)(A) - \zeta^r(0)(A + (\bar{S}^{r,m}(t), t)) \to \zeta(t)(A) - \zeta(0)(A + (tz^{-1}, t)).$$

Let $D_\vartheta(\mathscr{A}) = \{D \in \mathscr{A} : \vartheta(\partial_D) \neq 0\}$. Note that $D_\vartheta(\mathscr{A})$ is countable because $\vartheta(\cdot \times \mathbb{R})$ and $\vartheta(\mathbb{R}_+ \times \cdot)$ are probability measures. Since the function $s \mapsto S_{s,t}$ is strictly decreasing in $s$,

$$D_\vartheta(S) = \{s \in [0,t] : A + (S_{s,t} \pm 2\delta, t - s) \in D_\vartheta(\mathscr{A})\}$$

is also countable. For each integer $N > 1$, let $0 = t_0^N < t_1^N < \cdots < t_N^N = t$ be a partition of $[0,t]$ such that $t_j^N \notin D_\vartheta(S)$ for all $j = 1, \ldots, N-1$ and such that $\max_{j \leq N-1}(t_{j+1}^N - t_j^N) \to 0$ as $N \to \infty$. Let $I^r$ denote the second right-hand term in (5.5). Then

$$I^r = \sum_{j=0}^{N-1} \frac{1}{r} \sum_{i=r\bar{E}^{r,m}(t_j^N)+1}^{r\bar{E}^{r,m}(t_{j+1}^N)} 1_A^+(\bar{v}_i^{r,m}(t), \bar{l}_i^{r,m}(t)).$$

Suppose that $t_j^N \leq U_i^r r^{-1} - m \leq t_{j+1}^N$ for some $r > r_\delta$, some $j \leq N-1$ and some $i \in \{r\bar{E}^{r,m}(0)+1, \ldots, r\bar{E}^{r,m}(t)\}$. Then, by (5.8),

(5.12)            $$S_{t_{j+1}^N,t} - \delta \leq \bar{S}_{U_i^r r^{-1} - m,t}^{r,m} \leq S_{t_j^N,t} + \delta.$$

By the definitions (2.14), (2.15), (2.5) and (2.3),

$$(\bar{v}_i^{r,m}(t), \bar{l}_i^{r,m}(t)) = (v_i^r - \bar{S}_{U_i^r r^{-1} - m,t}^{r,m}, l_i^r r^{-1} - (t - (U_i^r r^{-1} - m))).$$

So, for $r > r_\delta$, (5.12) and the inequalities $1_A(\cdot - \delta, \cdot) \leq 1_A^+(\cdot, \cdot) \leq 1_A(\cdot + \delta, \cdot)$ yield

$$1_A^+(\bar{v}_i^{r,m}(t), \bar{l}_i^{r,m}(t)) \geq 1_A(v_i^r - (S_{t_j^N,t} + 2\delta), l_i^r r^{-1} - (t - t_j^N));$$

$$1_A^+(\bar{v}_i^{r,m}(t), \bar{l}_i^{r,m}(t)) \leq 1_A(v_i^r - (S_{t_{j+1}^N,t} - 2\delta), l_i^r r^{-1} - (t - t_{j+1}^N)).$$

This yields, for $r > r_\delta$,

$$I^r \geq \sum_{j=0}^{N-1} \frac{1}{r} \sum_{i=r\bar{E}^{r,m}(t_j^N)+1}^{r\bar{E}^{r,m}(t_{j+1}^N)} 1_A(v_i^r - (S_{t_j^N,t} + 2\delta), l_i^r r^{-1} - (t - t_j^N));$$

(5.13)

$$I^r \leq \sum_{j=0}^{N-1} \frac{1}{r} \sum_{i=r\bar{E}^{r,m}(t_j^N)+1}^{r\bar{E}^{r,m}(t_{j+1}^N)} 1_A(v_i^r - (S_{t_{j+1}^N,t} - 2\delta), l_i^r r^{-1} - (t - t_{j+1}^N)).$$

Using (2.16) and (2.18), we can rewrite (5.13) as

$$I^r \geq \sum_{j=0}^{N-1} \bar{\mathcal{L}}_{t_j^N, t_{j+1}^N}^{r,m}(A + (S_{t_j^N,t} + 2\delta, t - t_j^N));$$

(5.14)

$$I^r \leq \sum_{j=0}^{N-1} \bar{\mathcal{L}}_{t_j^N, t_{j+1}^N}^{r,m}(A + (S_{t_{j+1}^N,t} - 2\delta, t - t_{j+1}^N)).$$

By Definitions 4.17 and 4.15, (5.14) implies that, for $r > r_\delta$,

$$I^r \geq \sum_{j=0}^{N-1} (\alpha^r(t_{j+1}^N - t_j^N)\vartheta^r(A + (S_{t_j^N, t} + 2\delta, t - t_j^N)) - n(L, r)^{-1}),$$

$$I^r \leq \sum_{j=0}^{N-1} (\alpha^r(t_{j+1}^N - t_j^N)\vartheta^r(A + (S_{t_{j+1}^N, t} - 2\delta, t - t_{j+1}^N)) + n(L, r)^{-1}).$$

By (4.80) and the fact that $t_j^N \notin D_\vartheta(S)$ for all $j = 1, \ldots, N - 1$,

$$(5.15) \quad \liminf_{r \to \infty} I^r \geq \alpha \sum_{j=0}^{N-1} (t_{j+1}^N - t_j^N)\vartheta(A + (S_{t_j^N, t} + 2\delta, t - t_j^N)),$$

$$\limsup_{r \to \infty} I^r \leq \alpha \sum_{j=0}^{N-1} (t_{j+1}^N - t_j^N)\vartheta(A + (S_{t_{j+1}^N, t} - 2\delta, t - t_{j+1}^N)).$$

For $s \in [0, t]$ such that $s \notin D_\vartheta(S)$, the bounded convergence theorem implies that as $N \to \infty$,

$$(5.16) \quad \sum_{j=0}^{N-1} 1_{[t_j^N, t_{j+1}^N)}(s)\vartheta(A + (S_{t_j^N, t} + 2\delta, t - t_j^N)) \to \vartheta(A + (S_{s, t} + 2\delta, t - s)),$$

$$\sum_{j=0}^{N-1} 1_{[t_j^N, t_{j+1}^N)}(s)\vartheta(A + (S_{t_{j+1}^N, t} - 2\delta, t - t_{j+1}^N)) \to \vartheta(A + (S_{s, t} - 2\delta, t - s)).$$

Thus, the convergence in (5.16) holds for almost every $s \in [0, t)$. Let $N \to \infty$ in (5.15) and conclude from (5.16) and the bounded convergence theorem that

$$(5.17) \quad \liminf_{r \to \infty} I^r \geq \alpha \int_0^t \vartheta(A + (S_{s, t} + 2\delta, t - s)) \, ds,$$

$$\limsup_{r \to \infty} I^r \leq \alpha \int_0^t \vartheta(A + (S_{s, t} - 2\delta, t - s)) \, ds.$$

Let $\delta \to 0$ in (5.17). Since $D_\vartheta(\mathscr{A})$ is countable, both integrands in (5.17) converge almost everywhere on $[0, t]$ to

$$\vartheta(A + (S_{s, t}, t - s)) = \vartheta(A + ((t - s)z^{-1}, t - s)).$$

Conclude, using a change of variables, that

$$\lim_{r \to \infty} I^r = \alpha \int_0^t \vartheta(A + (sz^{-1}, s)) \, ds.$$

Combining this with (5.11) proves (5.2) for fixed $t \in [0, L)$ and $A \in \mathscr{A}$. To extend to all $B \in \mathscr{H}$, let $\mathscr{A}'$ be the set of $B \in \mathscr{H}$ for which (5.2) holds. Observe that $\mathscr{A}'$ is a $\lambda$-system: $\mathbb{H}_+ \in \mathscr{A}'$ because $\mathbb{H}_+ \in \mathscr{A}$; if $\{B_n\} \subset \mathscr{A}'$ satisfies $B_n \uparrow B$, then $B \in \mathscr{A}'$; if $B_1 \subset B_2$ are elements of $\mathscr{A}'$, then $B_2 \setminus B_1 \in \mathscr{A}'$. Also, observe that $\mathscr{A}$ is a $\pi$-system: if $A_1, A_2 \in \mathscr{A}$, then $A_1 \cap A_2 \in \mathscr{A}$. Since $\mathscr{A} \subset \mathscr{A}'$ and the $\sigma$-algebra generated by $\mathscr{A}$ is equal to $\mathscr{H}$, it follows that $\mathscr{A}' = \mathscr{H}$ by the Dynkin $\pi\lambda$-theorem (see, e.g., [1]). □

5.2. *Uniform approximation.* The next lemma states that with asymptotically high probability, the overlapping sections $\{\bar{\mathscr{Z}}^{r,m}(\cdot) : m \leq \lfloor rT \rfloor\}$ of the process $\bar{\mathscr{Z}}^r(\cdot)$ are uniformly approximated on $[0, L)$ by elements of the set $\mathscr{D}_L$.

Assume (A) and let $T > 1$ and $\eta \in (0, 1)$. Let $\{\mathscr{D}_L : L > 1\}$ be the subsets of $\mathbf{D}([0, \infty), \mathbf{M})$ specified by Definition 5.1.

LEMMA 5.4. *For each $L > 1$ and each $\varepsilon > 0$, the set $\mathscr{D}_L$ satisfies*

$$(5.18) \qquad \liminf_{r \to \infty} \mathbf{P}^r \left( \max_{m \leq \lfloor rT \rfloor} \inf_{\zeta(\cdot) \in \mathscr{D}_L} \sup_{t \in [0, L)} \mathbf{d}[\bar{\mathscr{Z}}^{r,m}(t), \zeta(t)] \leq \varepsilon \right) \geq 1 - \eta.$$

PROOF. Fix $L > 1$. By Definition 4.19, $\bar{\mathscr{Z}}_\omega^{r,m}(\cdot) 1_{[0,L)}(\cdot) \in \mathscr{D}_L^r$ for all $m \leq \lfloor rT \rfloor$ if and only if $\omega \in \Omega_L^r$. So, by Lemma 4.18,

$$(5.19) \qquad \liminf_{r \to \infty} \mathbf{P}^r(\bar{\mathscr{Z}}^{r,m}(\cdot) 1_{[0,L)}(\cdot) \in \mathscr{D}_L^r \text{ for all } m \leq \lfloor rT \rfloor) \geq 1 - \eta.$$

Thus, it suffices to show that

$$(5.20) \qquad \lim_{r \to \infty} \sup_{\zeta^r(\cdot) \in \mathscr{D}_L^r} \inf_{\zeta(\cdot) \in \mathscr{D}_L} \sup_{t \in [0, L)} \mathbf{d}[\zeta^r(t), \zeta(t)] = 0.$$

Suppose that (5.20) does not hold. Then there exists $\varepsilon > 0$, a subsequence $\tilde{\mathcal{R}} \subset \mathcal{R}$ and a sequence $\{\zeta^{\tilde{r}}(\cdot) : \tilde{r} \in \tilde{\mathcal{R}}\}$ with $\zeta^{\tilde{r}}(\cdot) \in \mathscr{D}_L^r$ for each $\tilde{r}$ such that

$$(5.21) \qquad \inf_{\tilde{r} \in \tilde{\mathcal{R}}} \inf_{\zeta(\cdot) \in \mathscr{D}_L} \sup_{t \in [0, L)} \mathbf{d}[\zeta^{\tilde{r}}(t), \zeta(t)] \geq \varepsilon.$$

By Theorem 4.20 and Definition 5.1, there exists a further subsequence $\{\tilde{r}_j\} \subset \tilde{\mathcal{R}}$ and a $\zeta(\cdot) \in \mathscr{D}_L$ such that

$$(5.22) \qquad \zeta^{\tilde{r}_j}(\cdot) \xrightarrow{J_1} \zeta(\cdot) \qquad \text{as } j \to \infty.$$

The path $\zeta(\cdot)$ is continuous on $[0, L)$ by Theorem 5.2(iii). Therefore,

$$(5.23) \qquad \lim_{j \to \infty} \sup_{t \in [0, L)} \mathbf{d}[\zeta^{\tilde{r}_j}(t), \zeta(t)] = 0,$$

contradicting (5.21). □

5.3. *Convergence to steady state.* This section establishes that when $L$ and $t$ are sufficiently large, $\zeta(t)$ is arbitrarily close to a steady state measure in $\mathbf{M}_\vartheta$, uniformly for all $\zeta(\cdot) \in \mathscr{D}_L$. The proof of this assertion exploits (5.2) and the constant total mass function $\langle 1, \zeta(\cdot) \rangle$.

Assume (A) and let $T > 1$ and $\varepsilon, \eta \in (0, 1)$. For each $L > 1$, let $\mathscr{D}_L$ be the set of local fluid limits specified by Definition 5.1.

THEOREM 5.5. *There exists $L^* > 1$ such that for all $\zeta(\cdot) \in \mathscr{D}_{L^*}$,*

$$(5.24) \qquad \sup_{t \in [L^* - 1, L^*)} \mathbf{d}[\zeta(t), \Delta_\vartheta \langle 1, \zeta(t) \rangle] \le \varepsilon.$$

PROOF. Let $M$ and $\mathbf{K}$ be given by Theorem 5.2. Choose $K > 0$ such that

$$(5.25) \qquad \sup_{\xi \in \mathbf{K}} \xi([KM^{-1}, \infty) \times \mathbb{R}) \le \varepsilon$$

and such that

$$(5.26) \qquad \alpha \int_K^\infty \vartheta([sM^{-1}, \infty) \times \mathbb{R}) \, ds \le \varepsilon;$$

such a $K$ exists because $\mathbf{K}$ is compact and $\langle \chi, \vartheta \rangle < \infty$ by (2.21). Choose $L^*$ such that

$$(5.27) \qquad L^* - 1 > K.$$

Fix $\zeta(\cdot) \in \mathscr{D}_{L^*}$ and $t \in [L^* - 1, L^*)$. By Theorem 5.2(v), $\langle 1, \zeta(\cdot) \rangle$ is constant and so $z = \langle 1, \zeta(0) \rangle$ satisfies

$$(5.28) \qquad \Delta_\vartheta \langle 1, \zeta(t) \rangle = \Delta_\vartheta \langle 1, \zeta(0) \rangle = \vartheta_e^z.$$

By Theorem 5.2(iv),

$$(5.29) \qquad z \le M.$$

Let $B \subset \mathbb{H}_+$ be a closed Borel set. By (5.28) and the definition of the Prohorov metric $\mathbf{d}[\cdot, \cdot]$, it suffices to verify the two inequalities

$$(5.30) \qquad \zeta(t)(B) \le \vartheta_e^z(B^\varepsilon) + \varepsilon,$$

$$(5.31) \qquad \vartheta_e^z(B) \le \zeta(t)(B^\varepsilon) + \varepsilon.$$

By Theorem 5.2(vi),

$$(5.32) \qquad \zeta(t)(B) = \zeta(0)(B + (tz^{-1}, t)) + \alpha \int_0^t \vartheta(B + (sz^{-1}, s)) \, ds.$$

Since $t \ge L^* - 1 > K$, (5.29) and (5.25) imply that

$$\zeta(0)(B + (tz^{-1}, t)) \le \zeta(0)([KM^{-1}, \infty) \times \mathbb{R})$$

$$\le \varepsilon.$$

Substitute this into (5.32); for the second term, enlarge $B$ to $B^\varepsilon$ and enlarge the domain of integration to conclude that

$$\zeta(t)(B) \leq \varepsilon + \alpha \int_0^\infty \vartheta(B^\varepsilon + (sz^{-1}, s))\, ds = \varepsilon + \vartheta_e^z(B^\varepsilon).$$

This proves (5.30). To prove (5.31), use Definition 2.1 to write

$$\vartheta_e^z(B) = \alpha \int_0^t \vartheta(B + (sz^{-1}, s))\, ds + \alpha \int_t^\infty \vartheta(B + (sz^{-1}, s))\, ds.$$

Enlarge $B$ to $B^\varepsilon$ in the first right-hand term; use $t \geq L^* - 1 > K$ and $z \leq M$ in the second right-hand term to obtain

$$\vartheta_e^z(B) \leq \alpha \int_0^t \vartheta(B^\varepsilon + (sz^{-1}, s))\, ds + \alpha \int_K^\infty \vartheta\big([sM^{-1}, \infty) \times \mathbb{R}\big)\, ds.$$

By (5.26),

$$\vartheta_e^z(B) \leq \alpha \int_0^t \vartheta(B^\varepsilon + (sz^{-1}, s))\, ds + \varepsilon.$$

Add an extra term and use Theorem 5.2(vi) to conclude that

$$\vartheta_e^z(B) \leq \zeta(0)(B^\varepsilon + (tz^{-1}, t)) + \alpha \int_0^t \vartheta(B^\varepsilon + (sz^{-1}, s))\, ds + \varepsilon$$

$$= \zeta(t)(B^\varepsilon) + \varepsilon,$$

which proves (5.31).  □

**6. Diffusion limit.**   The work of the previous sections is now combined to prove state space collapse, which immediately leads to the proof of Theorem 2.2.

Assume (A) and let $T > 1$.

THEOREM 6.1.   *As $r \to \infty$,*

$$(6.1) \qquad \sup_{t \in [0,T]} \mathbf{d}[\hat{\mathcal{Z}}^r(t), \Delta_\vartheta\langle 1, \hat{\mathcal{Z}}^r(t)\rangle] \Longrightarrow 0.$$

PROOF.   Let $\varepsilon, \eta \in (0,1)$. It suffices to show that

$$(6.2) \qquad \liminf_{r \to \infty} \mathbf{P}^r\left(\sup_{t \in [0,T]} \mathbf{d}[\hat{\mathcal{Z}}^r(t), \Delta_\vartheta\langle 1, \hat{\mathcal{Z}}^r(t)\rangle] \leq \varepsilon\right) \geq 1 - \eta.$$

By Theorem 5.5, there exists $L^* > 1$ such that the set $\mathscr{D}_{L^*}$ defined in Definition 5.1 satisfies

$$(6.3) \qquad \sup_{\zeta(\cdot) \in \mathscr{D}_{L^*}} \sup_{t \in [L^*-1, L^*)} \mathbf{d}[\zeta(t), \Delta_\vartheta\langle 1, \zeta(t)\rangle] \leq \frac{\varepsilon}{3}.$$

By Theorem 5.2(iv), there exists a positive $M$ such that for all $\zeta(\cdot) \in \mathscr{D}_{L^*}$,

$$(6.4) \qquad \sup_{t \in [0, L^*)} \langle 1, \zeta(t) \rangle \leq M.$$

Choose a positive $\delta \leq \varepsilon/36$ such that for each $z \in [0, M+1]$,

$$(6.5) \qquad \sup_{x \in \mathbb{R}_+} \langle 1_{[x, x+9\delta] \times \mathbb{R}}, \vartheta_e^z \rangle \leq \sup_{x \in \mathbb{R}_+} \langle 1_{[x, x+9\delta] \times \mathbb{R}}, \vartheta_e^{M+1} \rangle \leq \frac{\varepsilon}{6}$$

(see Definition 2.1). The map $\Delta_\vartheta : \mathbb{R}_+ \to \mathbf{M}_\vartheta$ is uniformly continuous on $[0, M+1]$ (Lemma 4.5). Thus, there exists a positive $\delta_1 \leq \delta$ such that

$$(6.6) \qquad \sup_{\substack{y, z \leq M+1, \\ |z-y| \leq \delta_1}} \mathbf{d}[\Delta_\vartheta z, \Delta_\vartheta y] \leq \delta.$$

For each $r \in \mathcal{R}$, define events

$$\Omega_1^r = \left\{ \max_{m \leq \lfloor rT \rfloor} \inf_{\zeta(\cdot) \in \mathscr{D}_{L^*}} \sup_{t \in [0, L^*)} \mathbf{d}[\bar{\mathcal{Z}}^{r,m}(t), \zeta(t)] < \delta_1 \right\},$$

$$\Omega_2^r = \{ \mathbf{d}[\bar{\mathcal{Z}}^{r,0}(0), \Delta_\vartheta \langle 1, \bar{\mathcal{Z}}^{r,0}(0) \rangle] < \delta_1 \},$$

$$\Omega_0^r = \Omega_1^r \cap \Omega_2^r.$$

By Lemma 5.4, (2.30) and (2.32),

$$(6.7) \qquad \liminf_{r \to \infty} \mathbf{P}^r(\Omega_0^r) \geq 1 - \eta.$$

Fix $\omega \in \Omega_0^r$ and assume henceforth in the proof that all random objects are evaluated at this $\omega$. By (6.7), it suffices to show that

$$\sup_{t \in [0, T]} \mathbf{d}[\hat{\mathcal{Z}}^r(t), \Delta_\vartheta \langle 1, \hat{\mathcal{Z}}^r(t) \rangle] \leq \varepsilon.$$

Note that if $t \in [L^*/r, T]$, then $\hat{\mathcal{Z}}^r(t) = \bar{\mathcal{Z}}^r(rt) = \bar{\mathcal{Z}}^{r,m}(s)$ for some $m \leq \lfloor rT \rfloor$ and some $s \in [L^*-1, L^*)$. If $t \in [0, L^*/r)$, then $\hat{\mathcal{Z}}^r(t) = \bar{\mathcal{Z}}^r(rt) = \bar{\mathcal{Z}}^{r,0}(s)$, for some $s \in [0, L^*)$. Thus, it suffices to show the two inequalities

$$(6.8) \qquad \max_{m \leq \lfloor rT \rfloor} \sup_{t \in [L^*-1, L^*)} \mathbf{d}[\bar{\mathcal{Z}}^{r,m}(t), \Delta_\vartheta \langle 1, \bar{\mathcal{Z}}^{r,m}(t) \rangle] \leq \varepsilon,$$

$$(6.9) \qquad \sup_{t \in [0, L^*)} \mathbf{d}[\bar{\mathcal{Z}}^{r,0}(t), \Delta_\vartheta \langle 1, \bar{\mathcal{Z}}^{r,0}(t) \rangle] \leq \varepsilon.$$

Fix $m \leq \lfloor rT \rfloor$ and $t \in [0, L^*)$. By the definition of $\Omega_1^r$, there exists $\zeta(\cdot) \in \mathscr{D}_{L^*}$ such that

$$(6.10) \qquad \mathbf{d}[\bar{\mathcal{Z}}^{r,m}(t), \zeta(t)] < \delta_1.$$

By definition of the Prohorov metric $\mathbf{d}[\cdot, \cdot]$, this implies that

$$|\langle 1, \zeta(t) \rangle - \langle 1, \bar{\mathcal{Z}}^{r,m}(t) \rangle| \leq \delta_1 \leq 1$$

and so by (6.4) and (6.6),

$$(6.11) \qquad \mathbf{d}[\Delta_\vartheta \langle 1, \zeta(t) \rangle, \Delta_\vartheta \langle 1, \bar{\mathcal{Z}}^{r,m}(t) \rangle] \leq \delta.$$

By (6.10) and (6.11), we have

$$
(6.12)
\begin{aligned}
\mathbf{d}[\bar{\mathcal{Z}}^{r,m}(t), \Delta_\vartheta \langle 1, \bar{\mathcal{Z}}^{r,m}(t) \rangle] &\leq \mathbf{d}[\bar{\mathcal{Z}}^{r,m}(t), \zeta(t)] \\
&\quad + \mathbf{d}[\zeta(t), \Delta_\vartheta \langle 1, \zeta(t) \rangle] \\
&\quad + \mathbf{d}[\Delta_\vartheta \langle 1, \zeta(t) \rangle, \Delta_\vartheta \langle 1, \bar{\mathcal{Z}}^{r,m}(t) \rangle] \\
&\leq 2\delta + \mathbf{d}[\zeta(t), \Delta_\vartheta \langle 1, \zeta(t) \rangle].
\end{aligned}
$$

If $t \in [L^* - 1, L^*)$, substituting (6.3) into (6.12) yields

$$\mathbf{d}[\bar{\mathcal{Z}}^{r,m}(t), \Delta_\vartheta \langle 1, \bar{\mathcal{Z}}^{r,m}(t) \rangle] \leq 2\delta + \frac{\varepsilon}{3} < \varepsilon,$$

which proves (6.8). To prove (6.9), it suffices to show that for $m = 0$ and $t \in [0, L^*)$, the $\zeta(\cdot)$ appearing in (6.10)–(6.12) also satisfies

$$(6.13) \qquad \mathbf{d}[\zeta(t), \Delta_\vartheta \langle 1, \zeta(t) \rangle] \leq \frac{\varepsilon}{3}.$$

By Theorem 5.2(v), $\langle 1, \zeta(\cdot) \rangle$ is constant. Let $z = \langle 1, \zeta(0) \rangle$ so that

$$(6.14) \qquad \Delta_\vartheta \langle 1, \zeta(t) \rangle = \Delta_\vartheta \langle 1, \zeta(0) \rangle = \vartheta_e^z.$$

Fix a closed Borel set $B \subset \mathbb{H}_+$. It suffices to show the two inequalities

$$(6.15) \qquad \zeta(t)(B) \leq \vartheta_e^z(B^{\varepsilon/3}) + \frac{\varepsilon}{3},$$

$$(6.16) \qquad \vartheta_e^z(B) \leq \zeta(t)(B^{\varepsilon/3}) + \frac{\varepsilon}{3}.$$

By (6.14), (6.10), the definition of $\Omega_2^r$ and (6.11),

$$
(6.17)
\begin{aligned}
\mathbf{d}[\zeta(0), \vartheta_e^z] &\leq \mathbf{d}[\zeta(0), \bar{\mathcal{Z}}^{r,0}(0)] \\
&\quad + \mathbf{d}[\bar{\mathcal{Z}}^{r,0}(0), \Delta_\vartheta \langle 1, \bar{\mathcal{Z}}^{r,0}(0) \rangle] \\
&\quad + \mathbf{d}[\Delta_\vartheta \langle 1, \bar{\mathcal{Z}}^{r,0}(0) \rangle, \vartheta_e^z] \\
&\leq \delta_1 + \delta_1 + \delta \\
&\leq 3\delta.
\end{aligned}
$$

By Theorem 5.2(vi), $\zeta(\cdot)$ satisfies (5.2). Apply (6.17) to the first term on the right of (5.2) to obtain

$$(6.18) \quad \zeta(t)(B) \leq \vartheta_e^z((B + (tz^{-1}, t))^{3\delta}) + 3\delta + \alpha \int_0^t \vartheta(B + (sz^{-1}, s)) \, ds.$$

Note that

$$(6.19) \qquad \vartheta_e^z((B + (tz^{-1}, t))^{3\delta}) \leq \vartheta_e^z(B^{3\delta} + (tz^{-1}, t))$$
$$+ \sup_{x \in \mathbb{R}_+} \langle 1_{[x, x+3\delta] \times \mathbb{R}}, \vartheta_e^z \rangle.$$

Apply (6.5) and use the relation $B \subset B^{3\delta} \subset B^{\varepsilon/3}$ to conclude from (6.18) and (6.19) that

$$(6.20) \quad \zeta(t)(B) \leq \vartheta_e^z(B^{\varepsilon/3} + (tz^{-1}, t)) + \alpha \int_0^t \vartheta(B^{\varepsilon/3} + (sz^{-1}, s)) \, ds + 3\delta + \frac{\varepsilon}{6}.$$

Inequality (6.15) follows from (6.20) and Definition 2.1. To show (6.16), write

$$\vartheta_e^z(B) = \vartheta_e^z(B + (tz^{-1}, t)) + \alpha \int_0^t \vartheta(B + (sz^{-1}, s)) \, ds.$$

Use (6.17) and enlarge $B$ to $B^{3\delta}$ to obtain

$$(6.21) \quad \vartheta_e^z(B) \leq \zeta(0)((B + (tz^{-1}, t))^{3\delta}) + 3\delta + \alpha \int_0^t \vartheta(B^{3\delta} + (sz^{-1}, s)) \, ds.$$

Note that

$$\zeta(0)((B + (tz^{-1}, t))^{3\delta}) \leq \zeta(0)(B^{3\delta} + (tz^{-1}, t)) + \sup_{x \in \mathbb{R}_+} \langle 1_{[x, x+3\delta] \times \mathbb{R}}, \zeta(0) \rangle$$

and so by (6.17),

$$(6.22) \qquad \zeta(0)((B + (tz^{-1}, t))^{3\delta}) \leq \zeta(0)(B^{3\delta} + (tz^{-1}, t))$$
$$+ \sup_{x \in \mathbb{R}_+} \langle 1_{[x, x+9\delta] \times \mathbb{R}}, \vartheta_e^z \rangle + 3\delta.$$

Apply (6.5) to deduce from (6.21) and (6.22) that

$$\vartheta_e^z(B) \leq \zeta(0)(B^{3\delta} + (tz^{-1}, t)) + \alpha \int_0^t \vartheta(B^{3\delta} + (sz^{-1}, s)) \, ds + 6\delta + \frac{\varepsilon}{6}.$$

Conclude from Theorem 5.2(vi) that

$$(6.23) \qquad \vartheta_e^z(B) \leq \zeta(t)(B^{3\delta}) + 6\delta + \frac{\varepsilon}{6},$$

which proves (6.16). This completes the proof of (6.13), which, together with (6.12), proves (6.9).  $\square$

PROOF OF THEOREM 2.2.   By Proposition 4.2,

$$\langle 1, \hat{Z}^r(\cdot) \rangle \Longrightarrow Z^*(\cdot) \qquad \text{as } r \to \infty.$$

Since the lifting map $\Delta_\vartheta$ is continuous (Lemma 4.5), the continuous mapping theorem implies that

$$\Delta_\vartheta \langle 1, \hat{\mathcal{Z}}^r(\cdot) \rangle \Longrightarrow \mathcal{Z}^*(\cdot) \qquad \text{as } r \to \infty.$$

Conclude from Theorem 6.1 and the "converging together lemma" ([2], Theorem 3.1) that

$$\hat{\mathcal{Z}}^r(\cdot) \Longrightarrow \mathcal{Z}^*(\cdot) \qquad \text{as } r \to \infty. \qquad \square$$

**Acknowledgment.** We are grateful to John Lehoczky and Steve Shreve for insightful discussions at the start of this work. We are grateful for the stimulating environment at EURANDOM, where part of this research was conducted. We also wish to thank two anonymous referees for their considerable effort and helpful suggestions.

DEPARTMENT OF MATHEMATICS
STANFORD UNIVERSITY
STANFORD, CALIFORNIA 94305-2125
USA
E-MAIL: gromoll@math.stanford.edu

DEPARTMENT OF MATHEMATICS
MARIA CURIE-SKŁODOWSKA UNIVERSITY
PL. MARII CURIE-SKŁODOWSKIEJ 1
20-031 LUBLIN
POLAND
E-MAIL: lkruk@hektor.umcs.lublin.pl